\documentclass[11pt]{preprint}
\usepackage[full]{textcomp}
\usepackage[osf]{newtxtext} 
\usepackage{colortbl}

\usepackage{upgreek}
\usepackage{comment}
\usepackage{amssymb}
\usepackage{lmodern}
\usepackage{mathtools}

\usepackage{hyperref}
\usepackage{breakurl}
\usepackage{aligned-overset}
\usepackage{mhenvs}
\usepackage{mhequ} 
\newcommand{\be}{\begin{equation*}}
\newcommand{\ee}{\end{equation*}}
\usepackage{mhsymb}
\usepackage{booktabs}
\usepackage{tcolorbox}
\usepackage[utf8]{inputenc}
\usepackage{longtable}
\usepackage{wrapfig}
\usepackage{subcaption}
\usepackage{epsfig}
\usepackage{microtype}
\usepackage[dvipsnames]{xcolor}
\usepackage{mathrsfs}
\usepackage{empheq}

\usepackage{centernot}
\usepackage{enumitem}
\usepackage{bm}
\usepackage{stackrel}
\usepackage{graphicx}
\usepackage{ulem}

\makeatletter
\newcommand{\globalcolor}[1]{%
  \color{#1}\global\let\default@color\current@color
}
\makeatother

\usetikzlibrary{calc}
\usetikzlibrary{decorations}
\usetikzlibrary{positioning}
\usetikzlibrary{shapes}
\usetikzlibrary{shapes.misc}

\tikzset{cross/.style={cross out, draw=black, fill=none, minimum size=2*(#1-\pgflinewidth), inner sep=0pt, outer sep=0pt}, cross/.default={2pt}}

\definecolor{blush}{rgb}{0.87, 0.36, 0.51}
	\definecolor{brightcerulean}{rgb}{0.11, 0.67, 0.84}
	\definecolor{greenryb}{rgb}{0.4, 0.69, 0.2}

\newif\ifdark
\darkfalse

\ifdark
\definecolor{darkred}{rgb}{0.9,0.2,0.2}
\definecolor{darkblue}{rgb}{0.7,0.3,1}
\definecolor{darkgreen}{rgb}{0.1,0.9,0.1}
\definecolor{franck}{rgb}{0,0.8,1}
\definecolor{pagebackground}{rgb}{.15,.21,.18}
\definecolor{pageforeground}{rgb}{.84,.84,.85}
\pagecolor{pagebackground}
\AtBeginDocument{\globalcolor{pageforeground}}
\definecolor{symbols}{rgb}{0,0.7,1}
\colorlet{connection}{red!80!black}
\colorlet{boxcolor}{blue!50}

\else

\definecolor{darkred}{rgb}{0.7,0.1,0.1}
\definecolor{darkblue}{rgb}{0.4,0.1,0.8}
\definecolor{darkgreen}{rgb}{0.1,0.7,0.1}
\definecolor{franck}{rgb}{0,0,1}
\definecolor{pagebackground}{rgb}{1,1,1}
\definecolor{pageforeground}{rgb}{0,0,0}
\colorlet{symbols}{blue!90!black}
\colorlet{connection}{red!30!black}
\colorlet{boxcolor}{blue!50!black}

\fi

\def\slash{\leavevmode\unskip\kern0.18em/\penalty\exhyphenpenalty\kern0.18em}
\def\dash{\leavevmode\unskip\kern0.18em--\penalty\exhyphenpenalty\kern0.18em}

\newcommand{\eqlaw}{\stackrel{\mbox{\tiny law}}{=}}

\DeclareMathAlphabet{\mathbbm}{U}{bbm}{m}{n}

\DeclareFontFamily{U}{BOONDOX-calo}{\skewchar\font=45 }
\DeclareFontShape{U}{BOONDOX-calo}{m}{n}{
  <-> s*[1.05] BOONDOX-r-calo}{}
\DeclareFontShape{U}{BOONDOX-calo}{b}{n}{
  <-> s*[1.05] BOONDOX-b-calo}{}
\DeclareMathAlphabet{\mcb}{U}{BOONDOX-calo}{m}{n}
\SetMathAlphabet{\mcb}{bold}{U}{BOONDOX-calo}{b}{n}

\setlist{noitemsep,topsep=4pt,leftmargin=1.5em}

\DeclareMathAlphabet{\BF}{OT1}{ptm}{bx}{n}
\DeclareMathAlphabet{\mathbbm}{U}{bbm}{m}{n}

\DeclareMathAlphabet{\mcb}{U}{BOONDOX-calo}{m}{n}
\SetMathAlphabet{\mcb}{bold}{U}{BOONDOX-calo}{b}{n}
\DeclareFontFamily{U}{mathx}{\hyphenchar\font45}
\DeclareFontShape{U}{mathx}{m}{n}{
      <5> <6> <7> <8> <9> <10>
      <10.95> <12> <14.4> <17.28> <20.74> <24.88>
      mathx10
      }{}
\DeclareSymbolFont{mathx}{U}{mathx}{m}{n}
\DeclareMathSymbol{\bigtimes}{1}{mathx}{"91}

\newcommand{\talloblong}{
\,\begin{tikzpicture}[baseline=-0.7ex]
    \draw[line width=.5pt] (0,-0.2) rectangle (0.08,0.20);
\end{tikzpicture}\,
}

\def\mcT{\mathcal{T}}

\tikzstyle{tinydots}=[dash pattern=on \pgflinewidth off \pgflinewidth]
\tikzstyle{superdense}=[dash pattern=on 4pt off 1pt]

\newcommand{\mcR}{\mathcal{R}}

\newcommand{\mcB}{\mathcal{B}}

\newcommand{\mcS}{\CS}

\newcommand{\mcL}{\mathcal{L}}

\newcommand{\mcP}{\mathcal{P}}

\newcommand{\e}{{\mathrm{e}}}

\newcommand{\beq}{\begin{equation}}
\newcommand{\eeq}{\end{equation}}

\newcommand{\mfc}{\mathfrak{c}}

\newcommand{\mfC}{\mathfrak{C}}

\newcommand{\mfM}{\mathfrak{M}}

\newcommand{\mfl}{\mathfrak{l}}

\newcommand{\mfs}{\mathfrak{s}}

\newcommand{\msB}{\mathscr{B}}

\newcommand{\msQ}{\mathscr{Q}}
\newcommand{\msP}{\mathscr{P}}
\newcommand{\msR}{\mathscr{R}}

\newcommand{\epmu}{{\eps,\mu}}
\newcommand{\epnu}{{\eps,\nu}}
\def\i{\mathrm{i}}
\def\x{\boldsymbol{x}}
\def\y{\boldsymbol{y}}
\def\z{\boldsymbol{z}}
\def\w{\boldsymbol{w}}

\def\n{\boldsymbol{n}}
\def\a{\boldsymbol{a}}
\def\b{\boldsymbol{b}}
\def\c{\boldsymbol{c}}

\newcommand{\p}{\boldsymbol{p}}

\newcommand{\Kappa}{\text{\Large$\kappa$}}

\DeclareSymbolFont{rmlargesymbols}{OMX}{mdbch}{m}{n}
\DeclareMathSymbol{\rmintop}{\mathop}{rmlargesymbols}{82}

\newcommand{\D}{\partial}
\newcommand{\dint}{\,\mathrm{d}}

\newcounter{theorem}

\newenvironment{DIFnomarkup}{}{} 

\theorembodyfont{\rmfamily}

\newfont{\indic}{bbmss12}

\makeatletter
\pgfdeclareshape{crosscircle}
{
  \inheritsavedanchors[from=circle] 
  \inheritanchorborder[from=circle]
  \inheritanchor[from=circle]{north}
  \inheritanchor[from=circle]{north west}
  \inheritanchor[from=circle]{north east}
  \inheritanchor[from=circle]{center}
  \inheritanchor[from=circle]{west}
  \inheritanchor[from=circle]{east}
  \inheritanchor[from=circle]{mid}
  \inheritanchor[from=circle]{mid west}
  \inheritanchor[from=circle]{mid east}
  \inheritanchor[from=circle]{base}
  \inheritanchor[from=circle]{base west}
  \inheritanchor[from=circle]{base east}
  \inheritanchor[from=circle]{south}
  \inheritanchor[from=circle]{south west}
  \inheritanchor[from=circle]{south east}
  \inheritbackgroundpath[from=circle]
  \foregroundpath{
    \centerpoint%
    \pgf@xc=\pgf@x%
    \pgf@yc=\pgf@y%
    \pgfutil@tempdima=\radius%
    \pgfmathsetlength{\pgf@xb}{\pgfkeysvalueof{/pgf/outer xsep}}%
    \pgfmathsetlength{\pgf@yb}{\pgfkeysvalueof{/pgf/outer ysep}}%
    \ifdim\pgf@xb<\pgf@yb%
      \advance\pgfutil@tempdima by-\pgf@yb%
    \else%
      \advance\pgfutil@tempdima by-\pgf@xb%
    \fi%
    \pgfpathmoveto{\pgfpointadd{\pgfqpoint{\pgf@xc}{\pgf@yc}}{\pgfqpoint{-0.707107\pgfutil@tempdima}{-0.707107\pgfutil@tempdima}}}
    \pgfpathlineto{\pgfpointadd{\pgfqpoint{\pgf@xc}{\pgf@yc}}{\pgfqpoint{0.707107\pgfutil@tempdima}{0.707107\pgfutil@tempdima}}}
    \pgfpathmoveto{\pgfpointadd{\pgfqpoint{\pgf@xc}{\pgf@yc}}{\pgfqpoint{-0.707107\pgfutil@tempdima}{0.707107\pgfutil@tempdima}}}
    \pgfpathlineto{\pgfpointadd{\pgfqpoint{\pgf@xc}{\pgf@yc}}{\pgfqpoint{0.707107\pgfutil@tempdima}{-0.707107\pgfutil@tempdima}}}
  }
}
\makeatother

\def\decorate#1#2{
        \ifnum#2>0
    		\foreach \count in {1,...,#2}{
	       	let
				\p1 = (sourcenode.center),
                \p2 = (sourcenode.east),
				\n1 = {\x2-\x1},
				\n2 = {1mm},
				\n3 = {(1.3+0.6*(\count-1))*\n1},
				\n4 = {0.7*\n1}
			in 
        		node[rectangle,fill=symbols,rotate=30,inner sep=0pt,minimum width=0.2*\n2,minimum height=\n2] at ($(sourcenode.center) + (\n3,\n4)$) {}
				}
		\fi
        \ifnum#1>0
    		\foreach \count in {1,...,#1}{
	       	let
				\p1 = (sourcenode.center),
                \p2 = (sourcenode.east),
				\n1 = {\x2-\x1},
				\n2 = {1mm},
				\n3 = {(1.3+0.6*(\count-1))*\n1},
				\n4 = {0.7*\n1}
			in 
        		node[rectangle,fill=symbols,rotate=-30,inner sep=0pt,minimum width=0.2*\n2,minimum height=\n2] at ($(sourcenode.center) + (-\n3,\n4)$) {}
				}
		\fi
}

\tikzset{
    dectriangle/.style 2 args={
        triangle,
        alias=sourcenode,
        append after command={\decorate{#1}{#2}}
    },
    dectriangle/.default={0}{0},
}

\tikzset{
	dot/.style={circle,fill=black,inner sep=0pt, minimum size=1mm},
	graydot/.style={circle,draw=gray,inner sep=0pt, minimum size=1mm},
	loopnode/.style={circle,draw=black,inner sep=0pt, minimum size=1.5mm},
	treenode/.style={circle,fill=black,inner sep=0pt, minimum size=1.5mm},
	charge/.style={circle,draw=black,inner sep=0pt, minimum size=3mm},
	loopline/.style={->,semithick,shorten >=1pt,shorten <=1pt},
	treeline/.style={semithick, densely dashed, shorten >=1pt,shorten <=1pt},
	homoge/.style={font=\scriptsize,draw=black,fill=white,inner sep=2pt},
	homo/.style={font=\scriptsize},
	KK/.style={thick,shorten >=1pt,shorten <=1pt, densely dashed, >=latex},
kerAlg/.style={thick},
kerAlg2/.style={semithick,dashed},
kerP/.style={very thick},
arrow/.style={->,thick,shorten >=3pt,shorten <=3pt},
renorm/.style={ultra thick},
	kernel/.style={semithick,shorten >=1pt,shorten <=1pt},
	}

\tikzset{
	cross/.style={path picture={ 
  		\draw[symbols]
			(path picture bounding box.south east) -- (path picture bounding box.north west) (path picture bounding box.south west) -- (path picture bounding box.north east);
		}},
root/.style={circle,fill=green!50!black,inner sep=0pt, minimum size=1.2mm},
        dot/.style={circle,fill=pageforeground,inner sep=0pt, minimum size=1mm},
        dotred/.style={circle,fill=pageforeground!50!pagebackground,inner sep=0pt, minimum size=2mm},
        var/.style={circle,fill=pageforeground!10!pagebackground,draw=pageforeground,inner sep=0pt, minimum size=3mm},
        kernel/.style={semithick,shorten >=2pt,shorten <=2pt},
        kernels/.style={snake=zigzag,shorten >=2pt,shorten <=2pt,segment amplitude=1pt,segment length=4pt,line before snake=2pt,line after snake=5pt,},
        rho/.style={densely dashed,semithick,shorten >=2pt,shorten <=2pt},
           testfcn/.style={dotted,semithick,shorten >=2pt,shorten <=2pt},
        renorm/.style={shape=circle,fill=pagebackground,inner sep=1pt},
        labl/.style={shape=rectangle,fill=pagebackground,inner sep=1pt},
        xic/.style={very thin,circle,draw=symbols,fill=symbols,inner sep=0pt,minimum size=1.2mm},
        g/.style={very thin,rectangle,draw=symbols,fill=symbols!10!pagebackground,inner sep=0pt,minimum width=2.5mm,minimum height=1.2mm},
        xi/.style={very thin,circle,draw=symbols,fill=symbols!10!pagebackground,inner sep=0pt,minimum size=1.2mm},
	xies/.style={very thin,rectangle,fill=green!50!black!25,draw=symbols,inner sep=0pt,minimum size=1.1mm},
	xiesf/.style={very thin,rectangle,fill=green!50!black,draw=symbols,inner sep=0pt,minimum size=1.1mm},
        xix/.style={very thin,crosscircle,fill=symbols!10!pagebackground,draw=symbols,inner sep=0pt,minimum size=1.2mm},
        X/.style={very thin,cross,rectangle,fill=pagebackground,draw=symbols,inner sep=0pt,minimum size=1.2mm},
	xib/.style={thin,circle,fill=symbols!10!pagebackground,draw=symbols,inner sep=0pt,minimum size=1.6mm},
	xie/.style={thin,circle,fill=green!50!black,draw=symbols,inner sep=0pt,minimum size=1.6mm},
	xid/.style={thin,circle,fill=symbols,draw=symbols,inner sep=0pt,minimum size=1.6mm},
	xibx/.style={thin,crosscircle,fill=symbols!10!pagebackground,draw=symbols,inner sep=0pt,minimum size=1.6mm},
	kernels2/.style={very thick,draw=connection,segment length=12pt},
	keps/.style={thin,draw=symbols,->},
	kepspr/.style={thick,draw=connection,->},
	krho/.style={thin,draw=symbols,superdense,->},
	krhopr/.style={thick,draw=connection,superdense},
	triangle/.style = { regular polygon, regular polygon sides=3},
	not/.style={thin,circle,draw=connection,fill=connection,inner sep=0pt,minimum size=0.5mm},
	diff/.style = {very thin,draw=symbols,triangle,fill=red!50!black,inner sep=0pt,minimum size=1.6mm},
	diff1/.style = {very thin,dectriangle={1}{0},fill=red!50!black,draw=symbols,inner sep=0pt,minimum size=1.6mm},
	diff2/.style = {very thin,dectriangle={1}{1},fill=red!50!black,draw=symbols,inner sep=0pt,minimum size=1.6mm},
		diffmini/.style = {very thin,rectangle,fill=black,draw=black,inner sep=0pt,minimum size=0.75mm},
	 kernelsmod/.style={very thick,draw=connection,segment length=12pt},
	 rec/.style = {very thin,rectangle,fill=black,draw=black,inner sep=0pt,minimum size=2mm},
	cerc/.style={very thin,circle,draw=black,fill=symbols,inner sep=0pt,minimum size=2mm},
	stars/.style={very thin,star,star points=6,star point ratio=0.5, draw=black,fill=red,inner sep=0pt,minimum size=0.7mm},
	>=stealth,
        }
        \tikzset{
root/.style={circle,fill=black!50,inner sep=0pt, minimum size=3mm},
        circ/.style={circle,fill=white,draw=black,very thin,inner sep=.5pt, minimum size=1.2mm},
        round1/.style={fill=white,outer sep = 0,inner sep=2pt,rounded corners=1mm,draw,text=black,thin,minimum size=1.2mm},
          circ1/.style={circle,fill=red!10,draw=red,very thin,inner sep=.5pt, minimum size=1.2mm},
        rect/.style={fill=white,outer sep = 0,inner sep=2pt,rectangle,draw,text=black,thin,minimum size=1.2mm},
        rect1/.style={fill=white,outer sep = 0,inner sep=2pt,rectangle,draw,text=black,thin,minimum size=1.2mm},
        round2/.style={fill=red!10,outer sep = 0,inner sep=2pt,rounded corners=1mm,draw,text=black,thin,minimum size=1.2mm},
       round3/.style={fill=blue!10,outer sep = 0,inner sep=2pt,rounded corners=1mm,draw,text=black,thin,minimum size=1.2mm}, 
        rect2/.style={fill=black!10,outer sep = 0,inner sep=2pt,rectangle,draw,text=black,thin,minimum size=1.2mm},
        dot/.style={circle,fill=black,inner sep=0pt, minimum size=1.2mm},
        dotred/.style={circle,fill=black!50,inner sep=0pt, minimum size=2mm},
        var/.style={circle,fill=black!10,draw=black,inner sep=0pt, minimum size=3mm},
        kernel/.style={semithick,shorten >=2pt,shorten <=2pt},
         diag/.style={thin,shorten >=4pt,shorten <=4pt},
        kernel1/.style={thick},
        kernels/.style={snake=zigzag,shorten >=2pt,shorten <=2pt,segment amplitude=1pt,segment length=4pt,line before snake=2pt,line after snake=5pt,},
		kernels1/.style={snake=zigzag,segment amplitude=0.5pt,segment length=2pt},
		rho1/.style={densely dotted,semithick},
        rho/.style={densely dashed,semithick,shorten >=2pt,shorten <=2pt},
           testfcn/.style={dotted,semithick,shorten >=2pt,shorten <=2pt},
           visible/.style={draw, circle, fill, inner sep=0.25ex},
        renorm/.style={shape=circle,fill=white,inner sep=1pt},
        labl/.style={shape=rectangle,fill=white,inner sep=1pt},
        xic/.style={very thin,circle,fill=symbols,draw=black,inner sep=0pt,minimum size=1.2mm},
        xi/.style={very thin,circle,fill=blue!10,draw=black,inner sep=0pt,minimum size=1.2mm},
	xib/.style={very thin,circle,fill=blue!10,draw=black,inner sep=0pt,minimum size=1.6mm},
	xie/.style={very thin,circle,fill=green!50!black,draw=black,inner sep=0pt,minimum size=1mm},
	xid/.style={very thin,circle,fill=symbols,draw=black,inner sep=0pt,minimum size=1.6mm},
	edgetype/.style={very thin,circle,draw=black,inner sep=0pt,minimum size=5mm},
	nodetype/.style={very thick,circle,draw=black,inner sep=0pt,minimum size=5mm},
	kernels2/.style={very thick,draw=connection,segment length=12pt},
clean/.style={thin,circle,fill=black,inner sep=0pt,minimum size=1mm},	not/.style={thin,circle,fill=symbols,draw=connection,fill=connection,inner sep=0pt,minimum size=0.8mm},
	>=stealth,
        }

\def\one{\BF{1}}
\def\eps{\varepsilon}
\newcommand{\snnorm}[1]{{\big\vert\kern-0.25ex\big\vert\kern-0.25ex\big\vert #1 
    \big\vert\kern-0.25ex\big\vert\kern-0.25ex\big\vert}}
\newcommand{\nnorm}[1]{{\Big\vert\kern-0.25ex\Big\vert\kern-0.25ex\Big\vert #1 
    \Big\vert\kern-0.25ex\Big\vert\kern-0.25ex\Big\vert}}
\newcommand{\bnnorm}[1]{{\Biggl\vert\kern-0.25ex\Biggl\vert\kern-0.25ex\Biggl\vert #1 
    \Biggl\vert\kern-0.25ex\Biggl\vert\kern-0.25ex\Biggl\vert}}

\DeclareMathAlphabet{\mathpzc}{OT1}{pzc}{m}{it}

\def\eqref#1{(\ref{#1})}

\makeatletter 
\newcommand*{\bigcdot}{}
\DeclareRobustCommand*{\bigcdot}{%
  \mathbin{\mathpalette\bigcdot@{}}%
}
\newcommand*{\bigcdot@scalefactor}{.5}
\newcommand*{\bigcdot@widthfactor}{1.15}
\newcommand*{\bigcdot@}[2]{%
  \sbox0{$#1\vcenter{}$}
  \sbox2{$#1\cdot\m@th$}%
  \hbox to \bigcdot@widthfactor\wd2{%
    \hfil
    \raise\ht0\hbox{%
      \scalebox{\bigcdot@scalefactor}{%
        \lower\ht0\hbox{$#1\bullet\m@th$}%
      }%
    }%
    \hfil
  }%
}
\makeatother

\tcbset
{colframe=boxcolor,colback=symbols!7!pagebackground,coltext=pageforeground,
fonttitle=\bfseries,nobeforeafter,center title,size=fbox,boxsep=1.5pt,
top=0mm,bottom=0mm,boxsep=0mm,tcbox raise base}

\def\mail#1{\burlalt{#1}{mailto:#1}}

\usepackage{thmtools}

\usepackage[margin=6em]{geometry}

\renewcommand{\leq}{\leqslant}
\renewcommand{\geq}{\geqslant}
\newcommand{\OH}{\mathrm{H}}
\newcommand{\St}{\mathrm{St}}
\newcommand{\highunderline}[2][3pt]{%
  \mathrel{\mathpalette\highunderlineaux{{#1}{#2}}}%
}
\newcommand{\highunderlineaux}[2]{%
  \highunderlineauxi#1#2%
}
\newcommand{\highunderlineauxi}[3]{%
  \begingroup
  \setbox0=\hbox{$#1#3$}%
  \rlap{\raisebox{#2}{\rule{\wd0}{0.4pt}}}%
  \box0
  \endgroup
}
\newcommand{\undcdot}{\highunderline[0.15pt]{\bigcdot}}
\newcommand{\ovcdot}{~\overline{\bigcdot}~}

\begin{document}

\title{The stochastic Ginzburg-Landau equation with harmonic potential in the full subcritical regime}

\author{Samaël Mackowiak$^{1,2}$ and Aurélien Minguella$^2$}
\institute{ 
 Univ Brest,  CNRS UMR 6205,  Laboratoire de Mathématiques de Bretagne Atlantique,  F-29200,  Brest,  France \and Université de Lorraine, CNRS, IECL, F-54000 Nancy, France
  \\
Email:\ \begin{minipage}[t]{\linewidth}
\mail{samael.mackowiak@univ-lorraine.fr}\\
\mail{aurelien.minguella@univ-lorraine.fr}
\end{minipage}}

\maketitle 

\begin{abstract}
    In this paper, we construct the complex Ginzburg-Landau dynamic associated with the harmonic oscillator, with a random force which is white in time and covers the full subcritical range in space. We are able to renormalise the non-linearity using the same constant as the one used to renormalise the problem associated to the Laplacian. By an extension of the problem to complex times, we are able to prove that the solution of the complex Ginzburg-Landau is equal in law to the solution of the purely dissipative case, \ie the associated complex $\Phi^4$ model, but along a complex axis of time.
\end{abstract}
\setcounter{tocdepth}{2}
\setcounter{secnumdepth}{4}
\tableofcontents

\section{Introduction}

We are interested in the local well-posedness of the stochastic complex Ginzburg-Landau equation with harmonic potential on $\R^d$
\begin{equation}\label{eq:bareGL}
    \D_tu+(\i+\gamma)(\OH u+\lambda|u|^{2}u)=\sqrt{2\gamma}\,\xi\,.
\end{equation}
With $\gamma>0$, $\lambda\in\R$, and $\OH=-\Delta+|x|^2$ is the harmonic oscillator. $\xi$ is a complex-valued mean-zero Gaussian noise with the covariance structure
\begin{equation}
    \E\big[\langle\xi,\phi\rangle\overline{\langle\xi,\psi\rangle}\big]=\int_{\R\times\R^d\times\R^d} \!\phi(t,x)\mathrm{Cov}(x-y)\overline{\psi(t,y)}\dint t\dint x\dint y\,,
\end{equation}
for any test functions $\phi,\psi\colon\Lambda\to\C$ and a covariance function $\mathrm{Cov}\colon\R^d\to\R$. We assume that $\mathrm{Cov}$ satisfies the estimate for $x\in\R^d$ and $\mfl\in\N^d$
\begin{equation}
    |\D_x^\mfl \mathrm{Cov}(x)|\lesssim\frac{\e^{-c|x|}}{|x|^{2(\beta-1)+|\mfl|}}\,,
\end{equation}
for some constant $c>0$. We choose a noise that is always white in time. This is essential to preserve some scaling properties of our system with respect to $\gamma$. With such a choice, the noise has space-time Hölder regularity $-\beta-\kappa$ for every $\kappa>0$. We impose that $\beta$ satisfies
\begin{equation}\label{eq:condbeta}
    2<\beta<3\wedge(2+d/4)\,.
\end{equation}
Let us explain this assumption. If $\beta\leq2$, the equation is not singular and can be treated via the Da Prato-Debussche trick \cite{DPD03,BDF22b}. The equation is scaling subcritical for $\beta<3$. Finally, the condition $\beta<2+d/4$ prevents a covariance blow-up. In the case where $\xi$ is a white noise \ie $\mathrm{Cov}=\delta$, $\beta=1+d/2$ and we find back the usual condition $2\leqslant d<4$. In that same case, the equation \eqref{eq:bareGL} has the $\Phi^4$ measure with harmonic potential and coupling parameter $\lambda$ for invariant measure. It formally reads
\begin{equation}\label{eq:baremes}
    \mathrm{d}\mu_{\lambda}(\phi)\propto\exp\left(-\int_{\R^d}\frac{1}{2}|\nabla\phi(x)|^2+\frac{1}{2}|x|^2|\phi(x)|^2+\frac{\lambda}{4}|\phi(x)|^4\dint x\right)\dint\phi\,,
\end{equation}
where $\mathrm{d}\phi$ is the \textit{formal} Lebesgue measure on functions $\phi\colon\R^d\to\C$. This measure is normalised to be a probability measure. It is a classical problem in the physics and SPDE literature that both the equation \eqref{eq:bareGL} and the measure \eqref{eq:baremes} do not make sense as so. A \textit{renormalisation} procedure has to be carried out, \ie for a mollification of the noise, $\xi_\eps$, we solve the equation
\begin{equation}
    \D_tu_\eps+(\i+\gamma)\big(\OH u_\eps+(\lambda|u_\eps|^{2}-\mfC_\eps^\lambda)u_\eps\big)=\sqrt{2\gamma}\,\xi_\eps\,,
\end{equation}
where $\mfC_\eps^\lambda$ is a renormalisation constant that does not depend on $\gamma$, and that diverges as $\eps$ goes to zero. We expect that with a good choice of this constant, $u_\eps$ converges towards a non-trivial and well-posed distribution.

\subsection{Main results}

We can now state the main results of this work. We consider the slightly more general equation
\begin{equation}\label{eq:solvedeqOH}
    \D_tu_\eps^\gamma+\gamma\OH u_\eps^\gamma=-\gamma(\lambda|u_\eps^\gamma|^{2}-\mfC^\lambda_\eps)u_\eps^\gamma+\sqrt{2\Re(\gamma)}\,\xi_\eps\,,
\end{equation}
where $\gamma$ is a complex number with positive real part. We then have the local well-posedness result.
\begin{theorem}\label{thm:maintheorem}
    Assume that $\beta$ satisfies \eqref{eq:condbeta}. Let $p,q\in[1,\infty]$, and take an initial condition $u_0\in\msB_{p,q}^{-\sigma}(\R^d)$, with $\sigma<\beta$. Then there exists a diverging deterministic constant $\mfC_\eps^\lambda$, unique up to a shift by a real number, and a random time $T_\gamma$ such that the solution $u_\eps^\gamma$ to \eqref{eq:solvedeqOH} converges in probability to some non-trivial distribution $u^\gamma\in L_t^\infty\msB_{p,q}^{-\beta+2}\big((0,T_\gamma)\times\R^d\big)$. This limit is independent of the choice of mollifier. Moreover, there is an extension of $u^1$ to complex times such that
    \begin{equation}\label{eq:eqlawOH}
        u^\gamma(t,x)\eqlaw u^1(\gamma t,x)\,.
    \end{equation}
\end{theorem}

\begin{proof}
    The local well-posedness is a consequence of the stochastic estimates of Proposition \ref{prop:estimatezeta}, as well as a fixed-point argument for the system \eqref{eq:flowsys}. The latter and the construction of the solution then follow Section 2.5.3 of \cite{FMM26}, slightly adapted to fit the present framework. The second fact is proven in Proposition \ref{prop:eqlawu}. The meaning that we give to the extension of $u^1$ to complex times is provided in Section \ref{sec:complex}.
\end{proof}

As a by-result of the proof of this first theorem, we also obtain a result of local well-posedness of the following equation on the torus $\T^d$, thus extending the result of \cite{HIN22} to the full subcritical regime.
\begin{equation}\label{eq:solvedeqtorus}
    \D_tu_\eps^\gamma-\gamma\Delta u_\eps^\gamma=-\gamma(|u_\eps^\gamma|^{2}-\tilde\mfC^\lambda_\eps)u_\eps^\gamma+\sqrt{2\Re(\gamma)}\,\xi_\eps\,,
\end{equation}

\begin{theorem}
    Assume that $\beta$ satisfies \eqref{eq:condbeta}. Let $p,q\in[1,\infty]$, and take an initial condition $u_0\in\mcB_{p,q}^{-\sigma}(\T^d)$, with $\sigma<\beta$. Then there exists a diverging deterministic constant $\tilde\mfC_\eps^\lambda$, unique up to a shift by a real number, and a random time $T_\gamma$ such that the solution $u_\eps^\gamma$ to \eqref{eq:solvedeqOH} converges in probability to some non-trivial distribution $u^\gamma\in L_t^\infty\mcB_{p,q}^{-\beta+2}\big((0,T_\gamma)\times\T^d\big)$. This limit is independent of the choice of mollifier. Moreover, there is an extension of $u^1$ to complex times such that
    \begin{equation}
        u^\gamma(t,x)\eqlaw u^1(\gamma t,x)\,.
    \end{equation}
\end{theorem}

\begin{proof}
    We do not prove the theorem as such, but the reader can follow the lines of the proof of Theorem \ref{thm:maintheorem}, and see that the present theorem is actually much easier to prove than the latter. 
\end{proof}

\begin{remark}
    A few remarks are in order.
    \begin{enumerate}
        \item The constants appearing in Theorem \ref{thm:maintheorem} are the same as those used to renormalise the following SPDE on $\R^d$.
        \begin{equation}
            \D_tu_\eps^\gamma-\gamma\Delta u_\eps^\gamma=-\gamma(|u_\eps^\gamma|^{2}-\mfC^\lambda_\eps)u_\eps^\gamma+\sqrt{2\Re(\gamma)}\,\xi_\eps\,.
        \end{equation}
        However, even if it should be true, we do not prove the well-posedness of this equation in the present paper.
        \item Note that in Theorem \ref{thm:maintheorem}, we use Besov spaces based on the harmonic oscillator. In negative regularity, they allow additional growth in large volume, avoiding the use of external weights.
    \end{enumerate}
    
\end{remark}

\subsection{Strategy of proof}

Our strategy is to adapt the tools developed on top of the flow approach of Duch \cite{Duc21} in \cite{FMM26} to the present context. In practice we will refer to many proofs in the latter, as they are essentially identical. However, we dealt with a linear equation therein, where the Ginzburg-Landau equation is not only non-linear, but also complex-valued. A first challenge is to find a well-suited ansatz for the effective force. A second one is to implement the complex conjugation in that ansatz, and in the associated hierarchy of flow equations, so that everything closes nicely. Regarding the renormalisation procedure, the main feature of \cite{FMM26} was to construct the stochastic terms that depend on the Mehler kernel (\ie the elementary solution to the linear equation driven by the harmonic oscillator) pathwise in the same terms but based on the elementary solution to the classical heat equation. This allowed the authors to show that the harmonic oscillator problem could be renormalised using the same renormalisation constants as the Laplacian problem. Technically speaking, the main observation is that at short times (and thus at small scales for the scale-decomposed kernel), the Mehler kernel is close to the heat kernel, so that the renormalisation, that is related only to small scales issues, should be the same in both cases. In practice, we use an adapted version of the following estimate
\begin{equation}
    \Big|M(t,x,y)-\e^{-t\frac{|x|^2+|y|^2}{4}}G(t,x-y)\Big|\lesssim t\big(1+|x|^2+|y|^2\big)\e^{-ct\frac{|x|^2+|y|^2}{4}}\,,
\end{equation}
for any $0<c<1$. The presence of the Gaussian weight on the left-hand side is necessary to get the Gaussian weight on the right-hand side, which is essential to tackle the infinite volume. What we can see there is that we gain some power of $t$ on the right-hand side, which will transfer into a gain of $2$ in power counting. The price to pay is additional quadratic growth in large volume. There is then a subtle interplay between these two phenomena, which was overcome in \cite{FMM26}.

In addition to the local well-posedness, our method allows us to obtain more interesting features on the SPDE \eqref{eq:solvedeqOH}. First, the fact that the renormalisation constant appearing therein factorises in $\gamma$ is highly non-trivial. If it is not so hard to show in the case of the equation in dimension two driven by a space-time white noise, proving it in full generality requires substantial developments. Moreover, we point out that it is essential to show this identity to preserve the invariance of the measure \eqref{eq:baremes} with the renormalisation. The central property that makes the proof of that fact possible is the following property of the Mehler kernel $M^\gamma$, the fundamental solution of $\D_t+\gamma\OH$,
\begin{equation}
    M^\gamma(t,x,y)=M^1(\gamma t,x,y)\,,
\end{equation}
and similarly for the heat kernel. One can see that utilising this fact to link the stochastic objects coming from the equation with parameter $\gamma$ to the ones coming from the equation with parameter $1$ leads us to work with complex times. As a by-product, we obtain the equality in law \eqref{eq:eqlawOH}. The meaning that we give to the right-hand side is quite specific. A priori, the stochastic complex Ginzburg-Landau equation does not have an immediate extension to complex times because of its lack of holomorphic structure. We explain in Section \ref{sec:complex} how we solve the equation along complex radii.

\subsection{Literature}

The complex Ginzburg-Landau equation \eqref{eq:bareGL} arises in the study of Bose-Einstein condensates as an effective model of condensation at non-zero temperature. When there is no dissipation ($\gamma=0$), \eqref{eq:bareGL} gives back the Gross-Pitaevskii equation with one-particle Hamiltonian given by $\OH$,
\begin{equation}\label{eq:GP}
    \i\D_tu-\OH u=\lambda|u|^{2}u\,.
\end{equation}
We refer to \cite{Gro61,DPS05} for the physical derivation of this model, and to \cite{BGM00,Pic10,Gol13,Pic15,BPS16,Rou20} and references therein for the mathematical derivation of the the Gross-Pitaevskii equation without potential, and to \cite{Chen13} for the harmonic oscillator case. Let us also mention the literature on the broader family of Hartree-type linearities \cite{KP10,Pic11}, 
\begin{equation}\label{eq:Hartree}
    \i\D_tu-\OH u=(w*|u|^{2})u\,,
\end{equation}
including \cite{dVJZ25,Mac26} which deal with Anderson operators. At non-zero temperature, the condensate interacts with non-condensated particles, which can be modelled by a random force $\xi$ \cite{DS01,GD03}. A common choice is to take $\xi$ as a space-time white noise. In that case, in order to prevent the energy from growing to infinity, a dissipative term $\gamma>0$ is added, and should imply convergence to the equilibrium, \ie the Gibbs measure. Remark that, in full generality, we should replace the prefactor $\sqrt{2\gamma}$ by $\sqrt{2\gamma k_B T}$, where $k_B$ is the Boltzmann constant and $T$ is the temperature of the condensate, and the Gibbs measure is then formally given by 
\begin{equation}\label{eq:baremestemperature}
    \mathrm{d}\mu_{\lambda,T}(\phi)\propto\exp\left(-\frac{1}{k_B T}\int_{\R^d}\frac{1}{2}|\nabla\phi(x)|^2+\frac{1}{2}|x|^2|\phi(x)|^2+\frac{\lambda}{4}|\phi(x)|^4\dint x\right)\dint\phi\,.
\end{equation}
For the sake of simplicity, we will assume $k_B T=1$, but the method would apply for any positive temperature. Remark that the expected equilibrium state is independent of the dissipation parameter $\gamma>0$, and is also expected to be invariant by the flow of \eqref{eq:GP}. It is also the same equilibrium state as the purely dissipative case, \ie the associated $\Phi^4$ model. The parameter $\lambda\in\R$ appearing in \eqref{eq:bareGL} is often called the coupling parameter. Its sign depends on the nature of the internal interaction of the condensate. When $\lambda\geqslant 0$, we say that the non-linearity is defocusing, or that the interaction is repulsive, and we expect \eqref{eq:bareGL} to be globally well-posed. When $\lambda< 0$, we say that the non-linearity is focusing, or that the interaction is attractive, and we cannot exclude \eqref{eq:bareGL} to have solutions blowing-up in finite time with positive probability. In the purely dispersive and focusing case ($\gamma=0$ and $\lambda<0$), \eqref{eq:GP} is mass-supercritical on $\R^3$, meaning there are solutions blowing up in finite time of any mass (see \cite{AH23} for example), and  mass-critical on $\R^2$, meaning there are solutions that blow up in finite time of any mass higher than the critical mass (see \cite{BHHZ23} for example).

The study of \eqref{eq:bareGL} goes back to Da Prato and Debussche \cite{DPD03} and the study of the periodic $\Phi^4_2$ model, where the authors constructed the renormalised flow and proved the invariance of the Gibbs measure. On the torus, Hoshino, Inahama, and Naganuma \cite{HIN22} and then Trenberth \cite{Tre19} and Hoshino \cite{Hos18} studied the complex Ginzburg-Landau equation with the Laplacian using regularity structures, proving the global well-posedness of the equation in dimensions $2$ and $3$. In presence of an harmonic potential, \eqref{eq:bareGL} was first studied by De Bouard, Debussche and Fukuizumi on $\R^d$ for $d\in\{1,2\}$. On $\R$, \eqref{eq:bareGL} is not singular and is proven to be globally well-posed on the support of the Gibbs measure as well as exponentially converging to equilibrium \cite{BDF18}. On $\R^2$, they obtain local well-posedeness of the defocusing Wick-ordered, \ie renormalised with a inhomogeneous counterterm, version of \eqref{eq:bareGL} as well as the existence of an invariant measure and of stationary martingale solutions \cite{BDF22a,BDF22b}. In the 3d case, only the purely dissipative case (\ie the harmonic $\Phi^4_3$ case) has been treated so far, in the work of Deya, Fukuizumi and Thomann, where the well-posedness of \eqref{eq:bareGL} is obtained through an inhomogeneous renormalisation procedure à la Mourrat-Weber \cite{DFT25}. They also obtained existence and, in some regimes, uniqueness of an invariant measure \cite{DFT26}. We also stress that, to the authors' knowledge, the present paper is the first one to renormalise \eqref{eq:bareGL} by constants.

When $\gamma$ goes to $0$, we expect the solution of \eqref{eq:bareGL} to converge to the solution of \eqref{eq:GP}, this is the so-called inviscid limit. So far, the literature appears to be only available for the Laplacian. In the deterministic case, the problem was studied, for example, by Wu \cite{Wu98}, Bebouche and Jüngle \cite{BJ00} or Huang and Wang \cite{HW08}. In the stochastic setting, the literature focuses essentially on the convergence of stationary solutions, in the sense of invariant measures, as the Gibbs measure \eqref{eq:baremestemperature} does not depend on the dissipation parameter $\gamma$. Kuksin and Shirikyan \cite{KS04} obtained a weak inviscid limit to a stationary process that is more regular than the Gibbs measure. More recently, Zine \cite{Zin22} obtained a strong inviscid limit with respect to the Gibbs measure on $\T^2$, and Nguyen \cite{Ngu26} obtained a long-time inviscid limit on $[0,1]$.

The flow approach has recently been introduced by Duch \cite{Duc21,Duc22} and has known a gain of interest in the last years. It follows the two major breakthroughs \cite{Hai14,GIP15} that have reshaped the field of stochastic PDEs. It takes its ideas from Wilsonian renormalisation group ideas that were already known to the field with the work of Kupiainen \cite{Kup15}. If the original paper tackled only equations with polynomial non-linearities, it has been extended to general non-linearities by Chandra and Ferdinand in \cite{CF24a} and later to rough paths \cite{CF24b}. See also a work by the second author with Bruned \cite{BM25} for an automated treatment of the ansatz and of the renormalisation with tree techniques. One of the main successes of the flow approach is the construction of the $\Phi^4$ measure in the full subcritical regime and in infinite volume by Duch, Gubinelli, and Rinaldi in \cite{DGR23}. Part of the framework that we use in the present work was established in that paper. Very recently, the approach has also been used by the authors, together with Ferdinand, to construct Anderson-type operators based on non-local operators, as well as the harmonic oscillator \cite{FMM26}. The techniques used in the present work to circumvent the difficulties arising from the harmonic oscillator were derived then.

\subsection{Notations}\label{sec:notation}

We will denote in a general manner that $a\lesssim b$ if there is a positive constant $C$, which does not depend on the parameter at stake, depending on the context, such that $a\leq Cb$. We set $\R_+$ the subset of $\R$ of non-negative numbers, and $\C_+$ the subset of $\C$ of numbers of positive real part with in addition $0$. We denote $\Lambda=(-\infty,1]\times\R^d$. Similarly, we denote $\Lambda_{T}=[0,T]\times\R^d$. We choose to use fractional parabolic scaling for the space-time in $\Lambda$. We set $\mfs=(2,1,\dots,1)$ and, for $z\in\SS$, $|z|_\mfs=|z_0|^\frac{1}{2}+\sum_{i=1}^d|z_i|$. In particular, the effective dimension is $|\mfs|=d+2$. In addition, we set for $k\in\N^{d+1}$, $|k|_\mfs=2 k_0+\sum_{i=1}^dk_i$. For $x\in\R^d$, we define the Japanese bracket as $\langle x\rangle=(1+|x|^2)^\frac{1}{2}$. For $\gamma\in\C_+$, we set $\gamma\cdot(t,x)=(\gamma t,x)$.\\

\textit{Weights.} We set for $\mu\in(0,1]$ two $\mu$-dependent Japanese brackets on $x\in\R^d$ and $z\in\SS$ by
\begin{equation}\label{mubracket}
\langle x\rangle_{\mu}=\left(1+\frac{|x|^2}{\mu^2}\right)^\frac{1}{2}~~\mathrm{and}~~\langle z\rangle_{\mu}=\left(1+\frac{|z|_\mfs^2}{\mu^2}\right)^\frac{1}{2}.
\end{equation}
For a collection of $n$ points $\z=(z_1,\dots,z_n)\in\Lambda^n$ we denote by $\St(\z)$ its Steiner diameter, \ie the length of the smallest tree connecting all these points, with respect to the parabolic distance $|\bigcdot|_\mfs$. We also set for such a collection
\begin{equation}
    \langle\z\rangle_\mu=\left(1+\frac{\St(\z)^2}{\mu^2}\right)^\frac{1}{2}\,.
\end{equation}
Finally we define some exponential weights that will be crucial in the study of the harmonic oscillator. For a collection $\x=(x_1,\dots,x_n)\in(\R^d)^n$ and $\mu\in(0,1]$, we set
\begin{equation}\label{eq:wieghtW}
    \mcb W_\mu(\x)=\exp\bigg(\frac{\mu^2}{8}\sum_{k=1}^n|x_k|^2\bigg)\,,~~\mathrm{and}~~\mcb W_\mu^{-1}(\x)=\exp\bigg(-\frac{\mu^2}{8}\sum_{k=1}^n|x_k|^2\bigg)\,.
\end{equation}

\textit{Combinatorial sets and multi-indices.} We define two sets of multi-indices to parametrize our ansatz
\begin{equation}
    \mfM=\{a=(a_\circ,\underline{a},\overline{a},a_\pm)\in\N^4\colon a_\circ\leq N_\star\,,a_\circ+\underline{a}+\overline{a}\in 3\N\,,~~ \underline{a}\leq a_\circ+2,~~ \overline{a}\leq a_\circ+1\}\,,
\end{equation}
with more precise sets $\mfM_+=\{a\in\mfM\colon a_\pm=-\underline{a}+\overline{a}+1\}$ and $\mfM_-=\{a\in\mfM\colon a_\pm=-\underline{a}+\overline{a}-1\}$, and set $\mfM_\pm=\mfM_+\cup\mfM_-$. $N_\star$ is defined in \eqref{eq:N0N1}. We define the extended set
\begin{equation}
    \mfM^\star=\{\a=(a_\circ,\underline{a},\overline{a},a_\pm,a_c,a^c)\in\N^6\colon a_\circ+\underline{a}+\overline{a}+a_c+a^c\in 3\N\,,~~ a_c+\underline{a}\leq a_\circ+2,~~a^c+\overline{a}\leq a_\circ+1\}\,,
\end{equation}
and we define $\mfM_+^\star$, $\mfM_-^\star$, and $\mfM_\pm^\star$ as $\mfM_+$, $\mfM_-$, and $\mfM_\pm$, but with $N_\star$ replaced by $N^\star$. We give some details about the combinatorial interpretation of the parametrisation of the ansatz:
\begin{itemize}
    \item $a_\circ$: Loop order of the term.
    \item $\underline{a}$: Number of legs without conjugation.
    \item $\overline{a}$: Number of legs with conjugation.
    \item $a_\pm$: Charge, \ie number of noises $\xi$ minus the number of conjugated noises $\overline{\xi}$.
    \item $a_c$: Number of inserted harmonic flows of the initial condition.
    \item $a^c$: Number of inserted conjugated harmonic flows of the initial condition.
\end{itemize}
For $a,b,c\in\mfM$, we define $\overline{a}=(a_\circ,\overline{a},\underline{a},-a_\pm)$ as well as 
\begin{equation}
    b+c=(b_\circ+c_\circ+1,\underline{b}+\underline{c}-1,\overline{b}+\overline{c},b_\pm+c_\pm)~~\mathrm{and}~~b-c=(b_\circ+c_\circ+1,\underline{b}+\underline{c},\overline{b}+\overline{c}-1,b_\pm-c_\pm)\,.
\end{equation}
Note that $a\in\mfM_+\Leftrightarrow\overline{a}\in\mfM_-$. For $\a,\b,\c\in\mfM^\star$, we define $\overline{\a}=(a_\circ,\overline{a},\underline{a},-a_\pm,a^c,a_c)$ as well as
\begin{equation}\label{eq:def+}
    \b+\c=(b_\circ+c_\circ+1,\underline{b}+\underline{c}-1,\overline{b}+\overline{c},b_\pm+c_\pm,b_c+c_c,b^c+c^c)\,,
\end{equation}
and
\begin{equation}\label{eq:def-}
    \b-\c=(b_\circ+c_\circ+1,\underline{b}+\underline{c},\overline{b}+\overline{c}-1,b_\pm-c_\pm,b_c+c^c,b^c+c_c)\,.
\end{equation}
Note that $\a\in\mfM_+^\star\Leftrightarrow\overline{\a}\in\mfM_-^\star$. We finally define the total number of legs $l(a)=\underline{a}+\overline{a}$ and the total number of initial conditions $c(\a)=a_c+a^c$.\\

\textit{Power counting.}
We will require a notion of degree that follows the rules of power counting. We let for $\a\in\mfM$
\begin{equation}
    \deg(a)=2(1-\beta/3)a_\circ+(\underline{a}+\overline{a})\beta/3-\beta \,,
\end{equation}
and for $\a\in\mfM^\star$,
\begin{equation}
    \deg(\a)=\deg(a_\circ,\underline{a},\overline{a},a_\pm)+c(a)(\beta/3-\sigma)\,.
\end{equation}
We also set
\begin{equs}\label{eq:N0N1}
    N^\star=\lfloor(\beta+3\sigma/2)/(3-\beta)\rfloor\,,~\text{and}~~ N_\star=\lfloor \beta/(3-\beta)\rfloor\,,
\end{equs}
which are the smallest integers such that, for $\a\in\mfM_+^\star$ and $a\in\mfM_+$,
\begin{equs}
    n>N^\star\Rightarrow \deg(\a)>0\,,~\mathrm{and}~~n>N_\star\Rightarrow \deg(a)>0\,.
\end{equs}

\textit{Kernels.} We define, for $(t,x)\in\SS$, the heat kernel
\begin{equation}\label{greenfunction}
    G(t,x)=\one_{t>0}\frac{1}{(4\pi t)^{d/2}}\exp\left(-\frac{|x|^2}{4t}\right)\,,
\end{equation}
the fundamental solution of the differential operator $\D_t-\Delta$, and the Mehler kernel
\begin{equation}\label{eq:Mehlerkernel}
    M(t,x,y)=\one_{t>0}\frac{1}{\big(2\pi\sinh(2t)\big)^{d/2}}\exp\left(-\frac{|x-y|^2}{4\tanh(t)}-\tanh(t)\frac{|x+y|^2}{4}\right)\,,
\end{equation}
which is the fundamental solution of $\D_t-\Delta+|x|^2$. For $\gamma\in\C_+$ we set $G^\gamma(t,x)=G(\gamma t,x)$ and $M^\gamma(t,x,y)=M(\gamma t,x,y)$ the fundamental solutions of $\D_t-\gamma\Delta$ and $\D_t+\gamma\OH$ respectively. We also define, for a scale parameter $\mu\in(0,1]$, and $(t,x)\in \R\times\R^d$,
\begin{equation}
    Q_\mu(x)=(1-\mu^2\Delta)^{-1}(x)\,,
\end{equation}
and
\begin{equation}
    K_{1,\mu}(t,x)=\mu^{-2} \e^{-\frac{t}{\mu^2}}\one_{(0,\infty)}(t)(1-\mu^2\Delta)^{-1}(x)\,,
\end{equation}
which are respectively the fundamental solutions of $(1-\mu^2\Delta)$ and $\mcP_\mu=(1+\mu^2\D_t)\otimes(1-\mu^2\Delta)$. For a given integer $N\in\N$, we set
\begin{equation}
K_{N,\mu}=K_{1,\mu}^{*N}.
\end{equation}
In practice, it will be sufficient to take $K_\mu=K_{N,\mu}$ with $N=5$. We also set $\mcR_\mu=\big(\mcP_\mu^{N}\big)^*$, where $*$ denotes the adjunction with respect to the $L^2_{t,x}(\Lambda)$ inner product. We also introduce some kernels associated with the harmonic oscillator.  Moreover, we define a family of kernels that are going to serve the same role as $K_\mu$ in the solution theory for the SPDE with the harmonic oscillator. First, let a family of positive numbers $a_n=1+n/(2(N^\star+1))$, for $0\leq n\leq N^\star+1$. $N^\star$ is defined in \eqref{eq:N0N1}. In particular $1\leq a_n\leq 3/2<2$. We let
\begin{equation}\label{eq:defJmu}
    J_{\mu,n}(t,x,y)=K_{\mu}(t,x-y)\e^{-\frac{a_n\mu^2|y|^2}{8}}\,,
\end{equation}
and $J_\mu=J_{\mu,N^\star+1}$. We also define $\mcT_{\mu,n}=(J_{\mu,n}^{-1})^*$. More precisely, it reads, for $f\colon\Lambda\to\R$,
\begin{equation}\label{eq:defTmu}
    (\mcT_{\mu,n}f)=\mcR_\mu\Big(\e^{\frac{a_n\mu^2|\bigcdot|^2}{8}}f\Big)\,.
\end{equation}
We set $\mcT_\mu=\mcT_{\mu,N^\star+1}$. Finally, we define $\msQ_\nu$, $\msP_\mu$, and $\msR_\mu$ the same way as $Q_\mu$, $\mcP_\mu$, and $\mcR_\mu$ but with the Laplacian $-\Delta$ replaced by the harmonic oscillator $-\Delta+|x|^2$.\\

\textit{Function spaces.} For $f\colon\R^d\to\C$ we define the norm for $p\in[1,\infty]$,
\begin{equation}\label{eq:defLp}
    \|f\|_{L^p_x}=\left(\int_{\R^d}|f(x)|^p\dint x\right)^\frac{1}{p},
\end{equation}
and $L^p_x(\R^d)$ the usual Lebesgue-based space associated to it. Similarly, for $f\colon\SS\to\C$, we define, for $p,q\in[1,\infty]$,
\begin{equation}
    \|f\|_{L_{t,x}^{p,q}}=\left(\int_\R\left(\int_{\R^d}|f(t,x)|^q\dint x\right)^\frac{p}{q}\dint t\right)^\frac{1}{p},
\end{equation}
and consequently the associated space. When $p=q$, we simply denote it $L_{t,x}^p(\SS)$ or $L_z^p(\SS)$.

For $E$ and $E'$ two Banach spaces, we denote $\mcL(E,E')$ the space of bounded linear maps from $E$ to $E'$. It is classically endowed with the operator norm $\|A\|_{\mcL(E,E')}=\sup_{\|x\|_E\leq1}\|Ax\|_{E'}$. We set $\mcL(E)=\mcL(E,E)$. We make this norm explicit for kernel operators from $L_{t,x}^{\infty,r}$ to $L_{t,x}^{\infty,p}$, whose action reads
\begin{equation}
    f\in L_{t,x}^{\infty,r}\mapsto\left((t,x)\mapsto\int_0^\infty\int_{\R^d}k(t-s,x,y)f(s,y)\dint y\dint s\right)\,,
\end{equation}
for some kernel $k\colon\R\times\R^d\times\R^d\to\R$. We will make the confusion between the operator and its kernel. In this context, we will denote $\|k\|_{\mcL(L_{t,x}^{\infty,r},L_{t,x}^{\infty,p})}$ the norm of the kernel operator. For $\eta<0$, $p\in[1,\infty]$, we define the Besov space $\mcB_{p,q}^{\eta}(\R^d)$ as the completion of compactly supported smooth functions on $\R^d$ with respect to the norm
\begin{equation}
\|\bigcdot\|_{\mcB_{p,q}^{\eta}}=\Big\|\mu^{-\eta}\|Q_{\mu}^{*\lceil-\eta\rceil}\bigcdot\|_{L_{x}^{p}}\Big\|_{L^q(\mathrm{d\mu/\mu})}\,.
\end{equation}
Moreover, for $\eta\in(0,1)$, we define the space Besov spaces $\mcB_{p,q}^\eta(\R^d)$ as the completion of compactly supported smooth functions on $\R^d$ with respect to the norm
\begin{equation}
    \|\bigcdot\|_{\mcB_{p,\infty}^{\eta}}=\|\bigcdot\|_{L^p_{x}}+\Big\|\mu^{-\eta}\|(Q_\mu-\mathrm{Id})\bigcdot\|_{L^p_{x}}\Big\|_{L^q(\mathrm{d\mu/\mu})}\,.
\end{equation}
For $\eta\in(0,1)$, $p\in[1,\infty]$, and an interval $I$, we define the Besov space $L_t^\infty\mcB_{p,\infty}^{\eta}(I\times\R^d)$ as the completion of compactly supported smooth functions on $I\times\R^d$ with respect to the norm
\begin{equation}
\|f\|_{L_t^\infty\mcB_{p,\infty}^{\eta}}=\|f\|_{L_{t,x}^{\infty,p}}+\sup_{t\in I}\sup_{\mu\in (0,1]}\mu^{-\eta}\|(Q_\mu-\mathrm{Id})f(t,\bigcdot)\|_{L_{x}^{p}}\,.
\end{equation}
We finally define some Besov spaces based on the harmonic oscillator. For $\eta<0$, they are given by the completion of $\mcS(\R^d)$ under the norm
\begin{equation}
    \|\bigcdot\|_{\msB_{p,q}^{\eta}}=\Big\|\mu^{-\eta}\|\msQ^{\lceil-\eta\rceil}_{\mu} \bigcdot\|_{L^p_x}\Big\|_{L^q(\mathrm{d\mu}/\mu)}\,,
\end{equation}
and for $0<\eta<1$, under the norm
\begin{equation}
    \|\bigcdot\|_{\msB_{p,q}^{\eta}}=\|\bigcdot\|_{L^p_x}+\Big\|\mu^{-\eta}\|(\msQ_\mu-\mathrm{Id}) \bigcdot\|_{L^p_x}\Big\|_{L^q(\mathrm{d\mu}/\mu)}\,.
\end{equation}
We define $L_t^\infty\msB_{p,\infty}^{\eta}(I\times\R^d)$ as $L_t^\infty\mcB_{p,\infty}^{\eta}(I\times\R^d)$, but with $Q_\mu$ replaced by $\msQ_\mu$.

\subsection*{Acknowledgments}

{\small
    S.M. was supported by the ANR project Smooth ANR-22-CE40-0017. A.M. gratefully acknowledges funding support from the European Research Council (ERC) through the ERC Starting Grant Low Regularity Dynamics via Decorated Trees (LoRDeT), grant agreement No.\ 101075208. Views and opinions expressed are however those of the author(s) only and do not necessarily reflect those of the European Union or the European Research Council. Neither the European Union nor the granting authority can be held responsible for them. 
} 

\section{The flow approach}

We gain some generality in our parameters by solving the slightly more general equation
\begin{equation}
    \D_tu+\gamma\OH u=-\gamma\lambda|u|^{2}u+\sqrt{2\Re(\gamma)}\,\xi\,,
\end{equation}
with $\gamma\in\C_+$. This slight change of notation will allow us to consider the equation with $\gamma=1$ as our reference case. We regularise the noise by convolving it with a mollifier only in space. We differ here a bit from what is usual in the literature but doing so will allow us to preserve some scaling properties of our SPDE. More precisely, for $\eps>0$ small, we set $\xi_\eps(t,x)=\eps^{-d}\int_{\R^d}\rho(\eps^{-1}(x-y))\xi(t,y)\dint y$, where $\rho$ is a smooth function compactly supported in the unit ball of $\R^d$ and that integrates to $1$. We then solve the bare equation
\begin{equation}\label{eq:bareeq}
    \D_tu_\eps^\gamma+\gamma\OH u_\eps^\gamma=-\gamma\lambda|u_\eps^\gamma|^{2}u_\eps^\gamma+\mfC_\eps^{\gamma,\lambda}u_\eps^\gamma+\sqrt{2\Re(\gamma)}\,\xi_\eps\,,
\end{equation}
with renormalisation \textit{constant}
\begin{equation}
    \mfC_\eps^{\gamma,\lambda}=\sum_{n=0}^{\lfloor \frac{N_\star+1}{3}\rfloor}\mfc_\eps^{3n-1,\gamma,\lambda}\,,
\end{equation}
and regularised initial data $u_0^\eps(x)=\eps^{-d}\int_{\R^d}\rho(\eps^{-1}(x-y))u_0(y)\dint y$. With the help of Proposition \ref{prop:scalingC}, we know more precisely the form of the renormalisation. Namely, it reads
\begin{equation}
    \D_tu_\eps^\gamma+\gamma\OH u_\eps^\gamma=-\gamma(\lambda|u_\eps^\gamma|^{2}-\mfC^\lambda_\eps)u_\eps^\gamma+\sqrt{2\Re(\gamma)}\,\xi_\eps\,,
\end{equation}
where $\mfC^\lambda_\eps=\mfC_\eps^{1,\lambda}$ is a diverging constant that does not depend on $\gamma$. One can also show the slightly refined statement $\mfc_\eps^{3n-1,\gamma,\lambda}=\gamma\lambda^n\mfc_\eps^{3n-1,1,1}$. We let the non-linear functional
\begin{equation}
    S_\eps^\gamma[u]=-\gamma\lambda|u|^{2}u+\mfC_\eps^{\gamma,\lambda}u+\sqrt{2\Re(\gamma)}\,\xi_\eps\,.
\end{equation}
We do not emphasize the dependence of this functional on $\lambda$ for the sake of conciseness. Letting $\psi_\eps(t,x)=\delta(t)u_0^\eps(x)$, we can rewrite \eqref{eq:bareeq} in mild form
\begin{equation}\label{eq:mildu}
    u_\eps^\gamma=M^\gamma(\one_{(0,\infty)}S_\eps^\gamma[u_\eps^\gamma]+\psi_\eps)\,.
\end{equation}
We recall that $M^\gamma$ is the Mehler kernel with parameter $\gamma$. We first need to tackle the initial condition. For that purpose, we let $v_\eps^\gamma=u_\eps^\gamma-M^\gamma\psi_\eps$. It satisfies the equation 
\begin{equation}
    v_\eps^\gamma=M^\gamma(\one_{(0,\infty)}S_\eps^\gamma[v_\eps^\gamma+M^\gamma\psi_\eps])\,.
\end{equation}
We define a translated force by
\begin{equation}
    T_\eps^\gamma[u]=S_\eps^\gamma[u+M^\gamma\psi_\eps]-\psi_\eps\,,
\end{equation}
as well as the non-stationary force by
\begin{equation}
    F_\eps^\gamma[v]=\one_{(0,\infty)}T_\eps^\gamma[v]\,.
\end{equation}
With these definitions, we can rewrite the SPDE \eqref{eq:mildu} under the form
\begin{equation}
    v_\eps^\gamma=M^\gamma F_\eps^\gamma[v_\eps^\gamma]\,.
\end{equation}
We crucially used that $\one_{(0,\infty)}\psi_\eps=0$. We now implement the effective scale and the associated cut-off Green function. Namely, we let, for $\mu\in[0,1]$, 
\begin{equation}
    M^{\gamma}_\mu(t,x) =  \chi_{\mu}^\gamma(t) M^{\gamma}(t,x)\,,
\end{equation} 
where $\chi_{\mu}^\gamma(t) = \chi(|t\gamma|/\mu^2)$ with the convention $ \chi_{0}^\gamma(t) = 1 $. Here, we have supposed that the time $t$ belongs to $[0,1]$. The map $\chi \colon \R_+ \rightarrow \R_+$ is a fixed smooth and non-decreasing function  such that $\chi |_{[0,1]} = 0$ and $\chi |_{[2,\infty)} = 1$. In particular, $M^{\gamma}_0=M^{\gamma}$. We also denote $\dot M^{\gamma}_\mu=\D_\mu M^{\gamma}_\mu$. 
One has the following key properties on the support of $M^\gamma_{\mu}$
\begin{equation}\label{cutscale}
	\mathrm{Supp}\,M^{\gamma}_{\mu}\subset[\mu^2/|\gamma|,1/|\gamma|]\times\R^d,\quad\mathrm{Supp}\,\dot M^{\gamma}_{\mu} \subset [\mu^2/|\gamma|,2\mu^2/|\gamma|] \times \R^d\,.
\end{equation}
Similarly, we define
\begin{equation}
    G^{\gamma}_\mu(t,x) =  \chi_{\mu}^\gamma(t) G^{\gamma}(t,x)\,,
\end{equation}
and we point out that such a kernel shares the same support properties as $M^\gamma$. Note that we have to put some dependence on $\gamma$ in $\chi_\mu^\gamma$ to preserve the following scaling properties of the kernels
\begin{equation}
    M_\mu^\gamma(t,x)=M_\mu(\gamma t,x),~~\mathrm{and}~~\dot M_\mu^\gamma(t,x)=\dot M_\mu(\gamma t,x)\,,
\end{equation}
and similarly for $G$. We can now define the effective solution
\begin{equation}
    v_\epmu^\gamma=M_\mu^\gamma F_\eps^\gamma[v_\eps^\gamma]\,.
\end{equation}
We define the effective force $F_\epmu^\gamma$ and a remainder $R_\epmu^\gamma$ by imposing that they satisfy the relation
\begin{equation}
    F_\eps^\gamma[v_\eps^\gamma]=F_\epmu^\gamma[v_\epmu^\gamma]+R_\epmu^\gamma\,,
\end{equation}
with initial condition $F_{\eps,0}^\gamma=F_\eps^\gamma$. It was shown in \cite{DGR23} that $(v_\epmu^\gamma,R_\epmu^\gamma)_{\mu\in[0,1]}$ satisfies then the system of equations
\begin{empheq}[left=\empheqlbrace]{alignat=2}\label{eq:flowsys}
    \begin{aligned}
    \D_\mu  v^{\gamma}_\epmu&=\dot M^{\gamma}_\mu\big(F^{\gamma}_{\eps,\mu}[v^{\gamma}_{\eps,\mu}]+R^{\gamma}_{\eps,\mu}\big)\,,\\
    \D_\mu  R^{\gamma}_\epmu&=- \mathrm{D}F^{\gamma}_\epmu[v_\epmu^\gamma] \dot M^{\gamma}_\mu R^{\gamma}_\epmu-(\D_\mu+\mathrm{D}F_\epmu^\gamma[v_\epmu^\gamma]\dot M_\mu^\gamma)F_\epmu^\gamma[v^{\gamma}_{\eps,\mu}]\,.
    \end{aligned}
\end{empheq}
Note that making the choice $R_\epmu^\gamma=0$ yields that $F_\epmu^\gamma$ satisfies the flow equation
\begin{equation}\label{eq:poleq}
    \D_\mu F_\epmu^\gamma+\mathrm{D}F_\epmu^\gamma\dot M_\mu^\gamma F_\epmu^\gamma=0\,.
\end{equation}
However, due to some choices of norms we have to make to tackle the infinite volume problem, it is a simpler path that we cannot take. However, we can see that solving the system \eqref{eq:flowsys} amounts to solving \eqref{eq:poleq} modulo an error. We solve the system by projecting it on the so-called \textit{effective force coefficients}. Namely, we make the ansatz
\begin{equation}
    F_\epmu^\gamma[\phi](z)=\sum_{\a\in \mfM^\star_+}\int_{\Lambda^{l(\a)}}\!\!\upzeta_\epmu^{\a,\gamma}(z;z_{1:\underline{a}}^+,z_{1:\overline{a}}^-)\prod_{j=1}^{\underline{a}}\phi(z_j^+)\prod_{j=1}^{\overline{a}}\overline{\phi}(z_j^-)\dint z^+_{1:\underline{a}}\dint z^-_{1:\overline{a}}\,.
\end{equation}
We have denoted $z_{1:\underline{a}}^+=(z_1^+,\dots,z_{\underline{a}}^+)$ and $z_{1:\overline{a}}^-=(z_1^-,\dots,z_{\overline{a}}^-)$. We will try as much as possible to drop the dependence of $z^+$ and $z^-$ on $\underline{a}$ and $\overline{a}$ to lighten the notations. The system \eqref{eq:flowsys} implies that these coefficients satisfy the hierarchy of flow equations
\begin{equs}\label{eq:hierarchyzeta}
    \D_\mu\upzeta_\epmu^{\a,\gamma}=-\sum_{\a=\b+\c} B_+\big(\dot{M}_\mu^\gamma,\upzeta_\epmu^{\b,\gamma},\upzeta_\epmu^{\c,\gamma}\big)-\sum_{\a=\b-\c}B_-\big(\dot{M}^{\overline{\gamma}}_\mu,\upzeta_\epmu^{\b,\gamma},\overline{\upzeta_\epmu^{\c,\gamma}}\big)\,,
\end{equs}
where $B_+$ and $B_-$ are defined below. We recall that the operations in the indexes were defined in \eqref{eq:def+} and \eqref{eq:def-}. Note that the system of flow equations above is hierarchical in the loop order $a_\circ$.
\begin{definition}
    Let $\uplambda^1\in \mcS'(\Lambda\times\Lambda^{l_+^1}\times\Lambda^{l_-^1})$, $\uplambda^2\in \mcS'(\Lambda\times\Lambda^{l_+^2}\times\Lambda^{l_-^2})$, and $G\in\mcS'(\Lambda^2)$. For $l_+=l_+^1+l_+^2-1$ and $l_-=l_-^1+l_-^2$, we define an element of $\mcS'(\Lambda\times\Lambda^{l_+}\times\Lambda^{l_-})$ by
    \begin{multline}
        B_+\big(G,\uplambda^1,\uplambda^2\big)(z;z_{1:l_+}^+,z_{1:l_-}^-)\\=\sum_{j=1}^{l_+^1}\int_{\Lambda^2}\uplambda^1(z;z_{1:l_+^1}^+,z_{1:l_-^1}^-) G(z_j^+,w)\uplambda^2(w;z_{l_+^1+1:l_+^1+l_+^2}^+,z_{l_-^1+1:l_-^1+l_-^2}^-)\dint z_j^+\dint w\,.
    \end{multline}
    Similarly, for $l_+=l_+^1+l_+^2$ and $l_-=l_-^1+l_-^2-1$, we define an element of $\mcS'(\Lambda\times\Lambda^{l_+}\times\Lambda^{l_-})$ by
    \begin{multline}
        B_-\big(G,\uplambda^1,\uplambda^2\big)(z;z_{1:l_+}^+,z_{1:l_-}^-)\\=\sum_{j=1}^{l_-^1}\int_{\Lambda^2}\uplambda^1(z;z_{1:l_+^1}^+,z_{1:l_-^1}^-) G(z_j^-,w)\uplambda^2(w;z_{l_+^1+1:l_+^1+l_+^2}^+,z_{l_-^1+1:l_-^1+l_-^2}^-)\dint z_j^-\dint w\,.
    \end{multline}
\end{definition}
While solving \eqref{eq:hierarchyzeta} is our ultimate goal, the effective force coefficients are not well-suited to study the renormalisation problem because, in view of the initial condition $ F^{\gamma}_\eps[v]=\one_{(0,\infty)}S^{\gamma}_\eps[v+M^{\gamma}\psi_\eps]$,
\begin{enumerate}
    \item The effective force coefficients are not stationary in time.
    \item The effective force coefficients depend on the initial condition.
    \item The effective force coefficients are not stationary in space.
\end{enumerate}
While the solution to 1. and 2. is now quite classical, solving 3. was the main novelty of \cite{FMM26}. We give more details about this construction. We define first $S_\epmu^\gamma$ ($S$ as stationary) as a possibly approximate solution to the flow equation driven by the heat kernel $G^\gamma$
\begin{equation}
    \D_\mu S_\epmu^\gamma+\mathrm{D}S_\epmu^\gamma\dot G_\mu^\gamma S_\epmu^\gamma=0\,,
\end{equation}
with initial condition $S_\eps^\gamma$. We then define $O_\epmu^\gamma$ ($O$ as oscillator) and $T_\epmu^\gamma$ ($T$ as translation) as possibly approximate solutions to the flow equation \eqref{eq:poleq}, with respective initial conditions $S_\eps^\gamma$ and $T_\eps^\gamma$. Note that then $S_{\eps,0}^\gamma=O_{\eps,0}^\gamma$. We expand these forces respectively over the \textit{stationary force coefficients} $(\upxi_{{\eps,\mu}}^{a,\gamma})$, the \textit{oscillator force coefficients} $(\upomega_\epmu^{a,\gamma})$, and the \textit{translated force coefficients} $(\uptau_{{\eps,\mu}}^{\a,\gamma})$, writing 
\begin{align}
    S^{\gamma}_{\eps,\mu}[\phi](z)&=\sum_{a\in\mfM_+}\int_{\Lambda^{l(a)}}\!\!\upxi_\epmu^{a,\gamma}(z;z_{1:\underline{a}}^+,z_{1:\overline{a}}^-)\prod_{j=1}^{\underline{a}}\phi(z_j^+)\prod_{j=1}^{\overline{a}}\overline{\phi}(z_j^-)\dint z^+_{1:\underline{a}}\dint z^-_{1:\overline{a}}\,,\\
    O^{\gamma}_{\eps,\mu}[v](z)&=\sum_{a\in\mfM_+}\int_{\Lambda^{l(a)}}\!\!\upomega_\epmu^{a,\gamma}(z;z_{1:\underline{a}}^+,z_{1:\overline{a}}^-)\prod_{j=1}^{\underline{a}}\phi(z_j^+)\prod_{j=1}^{\overline{a}}\overline{\phi}(z_j^-)\dint z^+_{1:\underline{a}}\dint z^-_{1:\overline{a}}\,,\\
    T^{\gamma}_{\eps,\mu}[v](z)&=\sum_{\a\in\mfM^\star_+}\int_{\Lambda^{l(\a)}}\!\!\uptau_\epmu^{\a,\gamma}(z;z_{1:\underline{a}}^+,z_{1:\overline{a}}^-)\prod_{j=1}^{\underline{a}}\phi(z_j^+)\prod_{j=1}^{\overline{a}}\overline{\phi}(z_j^-)\dint z^+_{1:\underline{a}}\dint z^-_{1:\overline{a}}\,.
\end{align}
Our construction follows this order.
\begin{equs}
    \upxi_\epmu^{a,\gamma}\to\upomega_\epmu^{a,\gamma}\to\uptau_\epmu^{\a,\gamma}\to\upzeta_\epmu^{\a,\gamma}\,.
\end{equs}
The construction of the stationary force coefficients is essentially probabilistic and relies on a cumulant analysis. See Section \ref{sec:proba}. The renormalisation is implemented at this step. As the coefficients $(\upxi_\epmu^{a,\gamma})$ are translation invariant, they are renormalised with \textit{constants}. On the opposite, the coefficients $(\upomega_\epmu^{a,\gamma})$ are not stationary in space since the Mehler kernel is not, but they are still stationary in time. We construct them pathwise and continuously from the coefficients $(\upxi_\epmu^{a,\gamma})$, thus preserving their renormalisation. More details are provided in Section \ref{sec:OH}. Finally, the last two steps are by now classical and we will not give much details about this part of the construction. We refer the reader to \cite[Section 2.1 and 2.2]{FMM26}.

The coefficients $\upomega_\epmu^{a,\gamma}$, $a=(a_\circ,\underline{a},\overline{a})\in\mfM_+$ are given by the hierarchy of flow equations
\begin{equation}\label{eq:hierarchyomega}
    \D_\mu\upomega_\epmu^{a,\gamma}=-\sum_{a=b+c} B_+\big(\dot{M}_\mu^\gamma,\upomega_\epmu^{b,\gamma},\upomega_\epmu^{c,\gamma}\big)-\sum_{a=b-c}B_-\big(\dot{M}^{\overline{\gamma}}_\mu,\upomega_\epmu^{b,\gamma},\overline{\upomega_\epmu^{c,\gamma}}\big)\,,
\end{equation}
and the boundary conditions
\begin{equation}
    \upomega^{(0,0,0,1),\gamma}_\epmu(z)=\sqrt{2\Re(\gamma)}\,\xi_\eps(z)\,,~\mathrm{and}~~\upomega_\epmu^{(0,2,1,0),\gamma}(z,z_1^+,z_2^+,z_1^-)=-\gamma\lambda\delta(z-z_1^+)\delta(z-z_2^+)\delta(z-z_1^-)\,,
\end{equation}
and, for $n\leq \lfloor \frac{N_\star+1}{3}\rfloor$
\begin{equation}
    \upomega_{\eps,0}^{(3n-1,1,0,0),\gamma}(z,z_1^+)=\mfc_\eps^{3n-1,\gamma,\lambda}\delta(z-z_1^+)\,.
\end{equation}
The coefficients $\upxi_\epmu^{a,\gamma}$ are defined in the same way but $M_\mu^\gamma$ is replaced by $G_\mu^\gamma$, namely
\begin{equation}\label{eq:hierarchyxi}
    \D_\mu\upxi_\epmu^{a,\gamma}=-\sum_{a=b+c} B_+\big(\dot{G}_\mu^\gamma,\upxi_\epmu^{b,\gamma},\upxi_\epmu^{c,\gamma}\big)-\sum_{a=b-c}B_-\big(\dot{G}^{\overline{\gamma}}_\mu,\upxi_\epmu^{b,\gamma},\overline{\upxi_\epmu^{c,\gamma}}\big)\,,
\end{equation}
with the same boundary conditions as $(\upomega_\epmu^{a,\gamma})$. In the above, we have used the convention $\dot G_\mu^\gamma(z,w)=\dot G_\mu^\gamma(z-w)$.

In order to state the analytic statements concerning these objects, we need to introduce further notation. We set , for a field $\uplambda_\mu(\z)=\uplambda_\mu(z;z_{1:l_+}^+,z_{1:l_-}^-)$,
\begin{equation}
    \big(K_\mu^{\otimes (l+1)}\uplambda_\mu\big)(w;w_{1:l_+}^+,w_{1:l_-}^-)=\int_{\Lambda^{l+1}}\!\!K_\mu(w-z)\prod_{\substack{\scriptscriptstyle1\leq i\leq l_+\\\scriptscriptstyle1\leq j\leq l_-}}\!\!K_\mu(w_i^+-z_i^+)K_\mu(w_j^--z_j^-)\uplambda_\mu(z;z_{1:l_+}^+,z_{1:l_-}^-)\dint\z\,.
\end{equation}
We have implicitly set $l=l_++l_-$. We will also use this notation for other kernels than $K_\mu$.
\begin{definition}
    We define for such a field the following norm
    \begin{equation}
        \|\uplambda_\mu\|_{\mu,p,\kappa,\varsigma}=\sup_{t\in\R}\left(\int_\Lambda\left|\langle w\rangle^{-\kappa}\int_{\Lambda^l}\langle \w\rangle^{\varsigma-\kappa}_\mu|\uplambda_\mu(\w)|\dint w^+\dint w^-\right|^p\dint x\right)^\frac{1}{p}\,,
    \end{equation}
    where we have denoted $w=(t,x)$. We recall that we use the notation $\w=(w;w^+,w^-)$. We further set the shorter notation
    \begin{equation}
        \|\bigcdot\|_{\mu,\kappa}=\|\bigcdot\|_{\mu,\infty,\kappa,\kappa},~~\|\bigcdot\|_\mu=\|\bigcdot\|_{\mu,\infty,0,0}\,,~\mathrm{and}~~\|\bigcdot\|_{\mu,p}=\|\bigcdot\|_{\mu,p,0,0}\,.
    \end{equation}
    Moreover, we measure the size of all the stationary force coefficients at once by setting
    \begin{equation}
        \talloblong \upxi^\gamma\talloblong_{\kappa,\varsigma}=\max_{a\in\mfM}\Big(\sup_{\epmu\in(0,1]}\mu^{-\deg(a)+\kappa}\|\upxi_\epmu^{a,\gamma}\|_{\mu,\infty,\kappa,\varsigma}\Big)\,.
    \end{equation}
\end{definition}

For the rest of the paper, we will fix $\kappa>0$ small enough and $\varsigma>0$ big enough. We define a parameter that increases through the induction $\kappa_a\eqdef 4^{a_\circ-N^\star-1}\kappa$, and we let $\kappa_0\eqdef4^{-N^\star-1}\kappa$. For conciseness, we also denote $\talloblong \upxi^\gamma\talloblong\eqdef\talloblong \upxi^\gamma\talloblong_{\kappa_0,\varsigma}$. We begin the results of the present section with the stochastic estimates on the stationary force coefficients.

\begin{proposition}\label{prop:estimatexi}
There is a choice of counterterms $(\mfc_\eps^{3n-1,\gamma,\lambda})_{n\leq \lfloor \frac{N_\star+1}{3}\rfloor}$ such that the following holds. For every $a\in\mfM_\pm$, uniform in $\epmu\in(0,1]$, we have
    \begin{equation}
        \|K_\mu^{\otimes(l(a)+1)}\upxi_\epmu^{a,\gamma}\|_{\mu,\kappa_a,\varsigma}\lesssim\talloblong\upxi^\gamma\talloblong\mu^{\deg(a)-\kappa_a}\,,
    \end{equation}
alongside the moment bound
    \begin{equation}
        \E\Big[\talloblong\upxi^\gamma\talloblong^k\Big]^\frac{1}{k}\lesssim 1\,,
    \end{equation}
for $k\geq1$. Moreover, we have convergence of the force coefficients in the above topology
    \begin{equation}
        \|K_\mu^{\otimes(l(a)+1)}(\upxi_\epmu^{a,\gamma}-\upxi^{a,\gamma}_{\eps',\mu})\|_{\mu,\kappa_a,\varsigma}\lesssim(\eps\vee\eps')^{\eta}\mu^{\deg(a)-\kappa_a-\eta}\,,
    \end{equation}
    for some small $\eta>0$.
\end{proposition}
\begin{proof}
    This is a consequence of the bound on cumulants of the stationary force coefficients \eqref{eq:boundcumu} and of a Kolmogorov-type argument that was developed in \cite{DGR23} and readapted in \cite{FMM26}. There is however a difference to notice. Indeed the complex-valued aspect of our equation changes things a bit. At the end of the proof of the Kolmogorov argument, one needs to consider the quantity
    \begin{equs}
        \E\left[\left|K_\mu^{\otimes(l(a)+1)}\upxi_\epmu^{a,\gamma}(w,w^+,w^-)\right|^2\right]=\E\left[K_\mu^{\otimes(l(a)+1)}\upxi_\epmu^{a,\gamma}(w,w^+,w^-)K_\mu^{\otimes(l(a)+1)}\overline{\upxi_\epmu^{a,\gamma}}(w,w^+,w^-)\right]\,.
    \end{equs}
    Then, thanks to Lemma \ref{lemma:conj}, we can write that the latter is equal to
    \begin{equs}
        \E\left[K_\mu^{\otimes(l(a)+1)}\upxi_\epmu^{a,\gamma}(w,w^+,w^-)K_\mu^{\otimes(l(a)+1)}\upxi_\epmu^{\overline{a},\overline{\gamma}}(w,w^-,w^+)\right]\,,
    \end{equs}
    and, after further operations, this is precisely of the form of the cumulants we have in Section \ref{sec:proba}. The rest of the argument follows then as usual. The second fact on the convergence follows the same lines of proof and we decide to not spend more time on it for conciseness.
\end{proof}

The construction is continued with the stochastic estimates on the oscillator force coefficients, that are deduced from the previous statement. We provide more details about this construction in Section \ref{sec:OH}.

\begin{proposition}\label{prop:estimateomega}
    For every $\kappa$ and $a\in\mfM_\pm$, we have that
    \begin{equation}
        \sup_{w\in\Lambda}\langle t\rangle^{-\kappa}\int_{\Lambda^{l(a)}}\langle \w\rangle_\mu^{1-\kappa}|\big(J_\mu^{\otimes(l(a)+1)}\upomega_\epmu^{a,\gamma}\big)(\w)|\dint w^+\dint w^-\lesssim\talloblong\upxi^\gamma\talloblong\mu^{\deg(a)-\kappa_a}\,,
    \end{equation}
    uniformly in $\eps$ and $\mu$. Moreover, as in Proposition \ref{prop:estimatexi}, a similar statement holds for $\upomega_\epmu^{a,\gamma}-\upomega_{\eps',\mu}^{a,\gamma}$ with a vanishing quantity on the right-hand side.
\end{proposition}

\begin{proof}
    This Proposition is proved using Lemma \ref{prop:estimateiota} and following the lines of the proof of Proposition 2.23 in \cite{FMM26}.
\end{proof}

Eventually, we may state the main result of the section, \ie the stochastic estimates on the effective force coefficients, that we will directly use to construct the remainder and then the solution. 

\begin{proposition}\label{prop:estimatezeta}
    For every $a\in\mfM_\pm$, the following holds uniformly in $\mu\in(0,1]$
    \begin{equation}
        \|J_\mu^{\otimes(l(a)+1)}\upzeta^{a,\gamma}_\epmu\|_\mu\lesssim\talloblong\upxi^\gamma\talloblong\mu^{\deg(a)-\kappa}\,,
    \end{equation}
    with the convention $\upzeta_\epmu^{a,\gamma}=\upzeta_\epmu^{(a_\circ,\underline{a},\overline{a},a_\pm,0,0),\gamma}$. Moreover, if $a_c\geq1$ or $a^c\geq1$, for $p\in[1,\infty]$,
    \begin{equation}
        \|J_\mu^{\otimes(l(\a)+1)}\upzeta_\epmu^{\a,\gamma}\|_{\mu,p}\lesssim\talloblong\upxi^\gamma\talloblong\|u_0^\eps\|_{\msB_{p,q}^{-\sigma}}^{c(\a)}\mu^{\deg(\a)-\kappa}\,.
    \end{equation}
    Moreover, similar statements hold for $\upzeta_\epmu^{\a,\gamma}-\upzeta_{\eps',\mu}^{\a,\gamma}$ with vanishing quantities on the right-hand side.
\end{proposition}

\begin{proof}
    The proof of this statement comes from Proposition \ref{prop:estimateomega}, and is essentially identical to the arguments developed in \cite[Section 2.5.2]{FMM26}, slightly modified to fit the ansatz given by the Ginzburg-Landau equation.
\end{proof}

\section{The harmonic oscillator problem}\label{sec:OH}

In this section, we show how to adapt the tools developed in \cite[Section 2.5]{FMM26} to our present context to construct the harmonic oscillator coefficients $(\upomega_\epmu^{a,\gamma})$ from the stationary force coefficients $(\upxi_\epmu^{a,\gamma})$. An interesting feature of our approach is that this construction is pathwise from the translation-invariant problem, on which we have already carried out a renormalisation procedure. We make the following ansatz on the harmonic oscillator coefficients
\begin{equation}\label{eq:ansatzomega}
    \upomega_\epmu^{a,\gamma}(\z)=\mcb W_\mu(\y)\upxi_\epmu^{a,\gamma}(\z)+\upiota_\epmu^{a,\gamma}(\z)\,.
\end{equation}
We have denoted by $\y$ the collection of the spatial components of $\z$. The idea behind that definition is that we can show that the coefficients $(\upomega_\epmu^{a,\gamma})$ can be decomposed as the sum of a term that is rough but translation invariant and a term that is smoother but not translation invariant. The main technical ingredient of this machinery is the observation that, at short times, the Mehler kernel and the heat kernel are very close, so we can hope that the difference is going to provide us a positive power of $t$ for small $t$, and thus of $\mu$. This estimate is formalized in Lemma \ref{lemma:diffMG}, and we gain $2$ of power counting, which is sufficient to make the coefficients $(\upiota_\epmu^{a,\gamma})$ irrelevant (\ie of positive degree). However, we have to pay the price of a quadratic growth at infinity that we need to absorb. We thus need to manage a subtle interplay between gain of regularity and growth penalisation. One can then observe in Lemma \ref{lemma:diffMG} that putting a Gaussian weight next to the heat kernel allows to get another Gaussian weight on the right-hand side. A way to make this appear is exactly the ansatz \eqref{eq:ansatzomega}.

The two hierarchies of flow equations \eqref{eq:hierarchyomega} and \eqref{eq:hierarchyxi} imply that $\upiota_\epmu^{a,\gamma}$ also satisfies a hierarchy of equations given by
\begin{equs}\label{eq:hierarchyiota}
    \D_\mu\upiota_\epmu^{a,\gamma}&=-\frac{\mu}{4}\bigg(|y|^2+\sum_{j=1}^{\underline{a}}|y_j^+|^2+\sum_{j=1}^{\overline{a}}|y_j^-|^2\bigg)\mcb W_\mu \upxi_\epmu^{a,\gamma}\\
    &\quad-\sum_{a=b+c}B_+\big(\dot M_\mu^\gamma-\mcb W_\mu \dot G_\mu^\gamma,\mcb W_\mu\upxi_\epmu^{b,\gamma},\mcb W_\mu\upxi_\epmu^{c,\gamma}\big)+B_+\big(\dot M_\mu^\gamma,\mcb W_\mu \upxi_\epmu^{b,\gamma},\upiota_\epmu^{c,\gamma}\big)\\
    &\qquad\qquad+B_+\big(\dot M_\mu^\gamma,\upiota_\epmu^{b,\gamma},\mcb W_\mu \upxi_\epmu^{c,\gamma}\big)+B_+\big(\dot M_\mu^\gamma,\upiota_\epmu^{b,\gamma},\upiota_\epmu^{c,\gamma}\big)\\
    &\quad-\sum_{a=b-c}B_-\big(\dot{M}^{\overline{\gamma}}_\mu-\mcb W_\mu \dot{G}^{\overline{\gamma}}_\mu,\mcb W_\mu\upxi_\epmu^{b,\gamma},\mcb W_\mu\upxi_\epmu^{\overline{c},\gamma}\big)+B_-\big(\dot{M}^{\overline{\gamma}}_\mu,\mcb W_\mu \upxi_\epmu^{b,\gamma},\upiota_\epmu^{\overline{c},\gamma}\big)\\
    &\qquad\qquad+B_-\big(\dot{M}^{\overline{\gamma}}_\mu,\upiota_\epmu^{b,\gamma},\mcb W_\mu \upxi_\epmu^{\overline{c},\gamma}\big)+B_-\big(\dot{M}^{\overline{\gamma}}_\mu,\upiota_\epmu^{b,\gamma}, \upiota_\epmu^{\overline{c},\gamma}\big)\,.
\end{equs}
We have denoted above $\mcb W_\mu \dot G_\mu^\gamma(z,w)=\mcb W_\mu(x,y)\dot G_\mu^\gamma(z-w)$. The induction is started by
\begin{equation}
    \upiota_\epmu^{(0,0,0,1),\gamma}(z)=\sqrt{2\Re(\gamma)}\Big(1-\e^{-\mu^2\frac{|y|^2}{8}}\Big)\xi_\eps(z)\,,
\end{equation}
as well as
\begin{equation}
    \upiota_\epmu^{(0,2,1,0),\gamma}(z,z_1^+,z_2^+,z_1^-)=-\gamma\Big(1-\e^{-\mu^2\frac{|y|^2}{2}}\Big)\delta(z-z_1^+)\delta(z-z_2^+)\delta(z-z_1^-)\,,
\end{equation}
and $\upiota_{\eps,0}^{a,\gamma}(\z)=0$ for any other $a\in\mfM_+$.

\begin{lemma}\label{prop:estimateiota}
     The following two statements hold
    \begin{equation}\label{eq:HR+2}
        \sup_{w\in\Lambda}\langle t\rangle^{-\kappa_a}\int_{\Lambda^{l(a)}}\langle x\rangle^{-2-\kappa_a}\langle \w\rangle^{1-\kappa_a}_\mu\big|\big(J_{\mu,a_\circ}^{\otimes(l(a)+1)}\upiota^{a,\gamma}_\epmu\big)(\w)\big|\dint w^+\dint w^-\lesssim  \talloblong\upxi^\gamma\talloblong\mu^{\deg(a)+2-\kappa_a}\,,
    \end{equation}
    and
    \begin{equation}\label{eq:HR+0}
        \sup_{w\in\Lambda}\langle t\rangle^{-\kappa_a}\int_{\Lambda^{l(a)}} \langle \w\rangle^{1-\kappa_a}_\mu\big|\big(J_{\mu,a_\circ+1}^{\otimes(l(a)+1)}\upiota^{a,\gamma}_\epmu\big)(\w)\big|\dint w^+\dint w^-\lesssim  \talloblong\upxi^\gamma\talloblong\mu^{\deg(a)-2\kappa_a}\,.
    \end{equation}
    We have used the general notation $w=(t,x)$.
\end{lemma}

\begin{proof}
    We let the reader convince themself that this proof can be adapted from that of \cite[Lemma 2.31]{FMM26}, but we recall the main ideas. The first observation is that the second inequality can be proven from the first one. For that, we can use the extra weight given in $J_{\mu,a_\circ+1}$ (compared to $J_{\mu,a_\circ}$) to convert the spatial weight in $x$ into $\mu^{-2-\kappa_a}$. From that point, the idea is then to prove \eqref{eq:HR+2} from both statements, using the hierarchy \eqref{eq:hierarchyiota}.

    We explain now how each of the terms appearing in this latter equation can bring an additional $+2$ of power counting, which would be enough to make $\upiota_\epmu^{a,\gamma}$ irrelevant. In the first term, we use for a second time the fact that we can sacrifice some extra weight to convert the quadratic term into $\mu^2$. This loss of weight is manageable since $J_{\mu,n}$ has a stronger Gaussian weight than $\mcb W_\mu$, and the stationary coefficients do not need additional decay. The second relies on the fact that Lemma \ref{lemma:diffMG} gives us more regularity, and the quadratic growth is absorbed by the quadratic weight appearing in $\eqref{eq:HR+2}$. The control of the two cross terms is based on on Lemma \ref{lemma:estimMmu}, and the fact that the terms in $\upiota$ give, by induction, $+2$ of power counting. Note that, at this step, the fact that the weights in $\langle \w\rangle_\mu$ have a decreasing exponent in the loop order is required to propagate the induction. This is due to the repartition of the weights. Eventually, even if the last term is expected to give $+4$ of power counting, we cannot directly use it since it would also give a quartic growth at infinity that we cannot control. The idea is then to use both \eqref{eq:HR+0} and \eqref{eq:HR+2} to get $+2$ of power counting and only quadratic growth at infinity, which is harmless.
\end{proof}

\section{Extension to complex times and scaling in $\gamma$}\label{sec:complex}

We extend our equation to complex times in the following way, assuming some holomorphic structure on the non-linearity
\begin{equation}\label{eq:eqholo}
    \D_{|\omega|} u(\omega,x)+\gamma\OH u(\omega,x)=-\gamma\lambda u(\omega,x)^2 \overline{u(\overline{\omega},x)} +\gamma \mfC_\eps^\lambda u(\omega,x)+\sqrt{2\Re(\gamma)}\,\xi_\eps(\omega,x)\,,
\end{equation}
where $\omega$ belongs to the sector
\begin{equation}
    \Sigma_{\gamma} =\left\{\omega\in\C_+,\; \gamma\omega,\overline{\gamma}\omega\in\C_+\right\}.
\end{equation}
The meaning that we give to \eqref{eq:eqholo} is by far not obvious, and we would like to give some details about it. First, we need to extend the noise and the Dirac delta to $\C_+\times\R^d$. We do it in the following way
\begin{equs}\label{eq:extensiondistribs}
    \xi(\omega,x)=\frac{\Re(\omega)}{|\omega|}\xi\big(\Re(\omega),x\big)\,,~\mathrm{and}~~\delta(\omega,x)=\frac{\Re(\omega)}{\omega}\delta\big(\Re(\omega),x\big)\,.
\end{equs}
This definition is natural with regard to the scaling properties of our SPDE. For example, with such a definition, for $(t,x)\in\R_+\times\R^d$, $\delta(\gamma t,x)=\gamma^{-1}\delta(t,x)$. With such a definition of the noise, one can see that $\sqrt{\Re(\gamma)}\,\xi_\eps(t,x)\eqlaw\gamma\xi_\eps(\gamma t,x)$. It turns out that this latter noise behaves much more nicely with respect to the complex structure we want to give to our equation. The solved equation then takes the form
\begin{equation}\label{eq:eqholochanged}
    \D_{|\omega|} u(\omega,x)+\gamma\OH u(\omega,x)=-\gamma\lambda u(\omega,x)^2 \overline{u(\overline{\omega},x)} +\gamma \mfC_\eps^\lambda u(\omega,x)+\sqrt{2}\gamma\xi_\eps(\gamma\omega,x)\,.
\end{equation}
Doing so, we solve an SPDE that is not our original one, but that will coincide for $\omega\in\R_+$, in the sense that they have the same law, which is enough for our purposes. More rigorously, the SPDE \eqref{eq:eqholochanged} is solved in the mild form 
\begin{equation}\label{eq:complexmildu}
    u_\eps^\gamma(\e^{\i\theta}\cdot z) = M^{\e^{\i\theta}\gamma}(\one_{(0,\infty)}\e^{\i\theta}\tilde S_\eps^\gamma[u_\eps^\gamma(\e^{\i\theta}\cdot\bigcdot)]+\psi_\eps)(z)\,,
\end{equation}
which extends \eqref{eq:mildu} to complex times along radii, where $\tilde S_\eps^\gamma$ is this time given by the right-hand side of \eqref{eq:eqholochanged} with $\omega=\e^{\i\theta}t$. We recall that $\cdot$ was defined in Section \ref{sec:notation}. We stress that, because of its structure, \eqref{eq:eqholochanged} must be solved with \eqref{eq:complexmildu} simultaneously in $\pm\theta$. This fact forces us to rather solve the system
\begin{empheq}[left=\empheqlbrace]{equation}\label{eq:eqholobis}
    \begin{aligned}
    \D_{|\omega|} u_+(\omega,x)+\gamma\OH u_+(\omega,x)&=-\gamma\lambda u_+(\omega,x)^2 u_-(\omega,x) +\gamma \mfC_\eps^\lambda u_+(\omega,x)+\sqrt{2}\gamma\xi_\eps(\gamma\omega,x)\,,\\
    \D_{|\omega|} u_-(\omega,x)+\overline{\gamma}\OH u_-(\omega,x)&=-\overline{\gamma}\lambda u_-(\omega,x)^2 u_+(\omega,x) +\overline{\gamma} \mfC_\eps^\lambda u_-(\omega,x)+\sqrt{2}\overline{\gamma}\overline{\xi_\eps(\overline{\gamma\omega},x)}\,,\\
    u_+(0)=\overline{u_-(0)}&=\psi_\eps\,,
    \end{aligned}
\end{empheq}
Then, if the mild solution $(u_+,u_-)$ of \eqref{eq:eqholobis} is unique, we immediately see that $u_-(\omega)=\overline{u_+(\overline{\omega})}$, \ie $u_+$ is a mild solution of \eqref{eq:eqholochanged}. Conversely, if $u_+$ is a mild solution to \eqref{eq:eqholochanged} and $u_-(\omega)=\overline{u_+(\overline{\omega})}$, then $(u_+,u_-)$ is a mild solution of \eqref{eq:eqholobis}. Note that our definition of the noise in complex times also gives the more precise equality in law $\big(\sqrt{\Re(\gamma)}\,\xi_\eps(t,x),\sqrt{\Re(\gamma)}\,\overline{\xi_\eps}(s,y)\big)\eqlaw\big(\gamma\xi_\eps(\gamma t,x),\overline\gamma\xi_\eps(\overline\gamma s,y)\big)$, which ensures the coherence of the extension of the noise with the system \eqref{eq:eqholobis}. Then, solving \eqref{eq:eqholobis} along two symmetric complex radii simultaneously gives in turn the new system
\begin{empheq}[left=\empheqlbrace]{alignat=2}\label{eq:eqholoray}
    \begin{aligned}
    \D_{t} u_++\gamma\e^{\i\theta}\OH u_+&=-\gamma\e^{\i\theta}\lambda u_+^2 u_- +\gamma \e^{\i\theta}\mfC_\eps^\lambda u_++\sqrt{2}\,\gamma\e^{\i\theta}\,\xi_\eps\big(\gamma\e^{\i\theta} \cdot \bigcdot\big)\,,\\
    \D_{t} u_-+\overline{\gamma\e^{\i\theta}}\OH u_-&=-\overline{\gamma\e^{\i\theta}}\lambda u_-^2 u_+ +\overline{\gamma\e^{\i\theta}} \mfC_\eps^\lambda u_-+\sqrt{2}\,\overline{\gamma\e^{\i\theta}}\,\overline{\xi_\eps\big(\overline{\gamma\e^{\i\theta}}\cdot \bigcdot\big)}\,,\\
    u_+(0)=\overline{u_-(0)}&=\psi_\eps\,,
    \end{aligned}
\end{empheq}
where this time the variable is $z=(t,x)$, with $t$ a true real variable. Note that remarkably, this system is the same as \eqref{eq:eqholobis}, but with $\gamma$ replaced by $\gamma \e^{\i\theta}$ and solved with real times. It is also coherent with Proposition \ref{prop:scalingC}, since $\mfC_\eps^{\gamma \e^{\i\theta},\lambda}=\gamma \e^{\i\theta}\mfC_\eps^{1,\lambda}$. We postpone the details of the resolution of \eqref{eq:eqholobis} to Section \ref{sec:system} in the appendix.

We begin with a lemma that will be key to prove the main result of this section, and that will bring many simplifications throughout the paper. 
\begin{lemma}\label{lemma:conj}
    For any $a\in\mfM_\pm$,
    \begin{equation}
        \overline{\upxi_\epmu^{a,\gamma}}(z;z^+,z^-)=\upxi_\epmu^{\overline{a},\overline{\gamma}}(z;z^-,z^+),~~\mathrm{and}~~\overline{\upomega_\epmu^{a,\gamma}}(z;z^+,z^-)=\upomega_\epmu^{\overline{a},\overline{\gamma}}(z;z^-,z^+)\,,
    \end{equation}
    and for any $a\in\mfM^\star_\pm$,
    \begin{equation}
        \overline{\uptau_\epmu^{\a,\gamma}}(z;z^+,z^-)=\uptau_\epmu^{\overline{\a},\overline{\gamma}}(z;z^-,z^+),~~\mathrm{and}~~\overline{\upzeta_\epmu^{\a,\gamma}}(z;z^+,z^-)=\upzeta_\epmu^{\overline{\a},\overline{\gamma}}(z;z^-,z^+)\,.
    \end{equation}
\end{lemma}

\begin{proof}
    The proof works by induction in every case. We start with $\upxi_\epmu^{a,\gamma}$. For the base case, we need to check two terms. We have first
    \begin{equs}
        \overline{\upxi^{(0,0,0,1),\gamma}_\epmu}(z)=\sqrt{2\Re(\gamma)}\,\overline{\xi_\eps(z)}=\upxi^{(0,0,0,-1),\overline{\gamma}}_\epmu(z)\,,
    \end{equs}
    and then
    \begin{equs}
        \overline{\upxi_\epmu^{(0,2,1,0),\gamma}}(z,z_1^+,z_2^+,z_1^-)=-\overline{\gamma}\lambda\delta(z-z_1^+)\delta(z-z_2^+)\delta(z-z_1^-)=\upxi_\epmu^{(0,1,2,0),\overline{\gamma}}(z;z_1^-,z_1^+,z_2^+)\,.
    \end{equs}
    We now turn to the induction part. For that we use the hierarchy \eqref{eq:hierarchyxi}. We write
    \begin{equs}
        \overline{\upxi_\epmu^{a,\gamma}}(z;z^+,z^-)=\overline{\upxi_{\eps,0}^{a,\gamma}}(z;z^+,z^-)+\int_0^\mu\D_\nu\overline{\upxi_\epnu^{a,\gamma}}(z;z^+,z^-)\dint\nu\,.
    \end{equs}
    The first term is non-zero only for $a=(3n-1,1,0,0)$, and it is a consequence of Proposition \ref{prop:scalingC} that for every $n$, $\overline{\mfc_\eps^{3n-1,\gamma,\lambda}}=\overline{\gamma}\mfc_\eps^{3n-1,1,\lambda}=\mfc_\eps^{3n-1,\overline{\gamma},\lambda}$. For the second term, one needs to notice that
    \begin{equs}
        \overline{B_+\big(\dot{G}_\nu^\gamma,\upxi_\epnu^{b,\gamma},\upxi_\epnu^{c,\gamma}\big)}(z;z^+,z^-)=B_-\big(\dot{G}_\nu^{\overline{\gamma}},\overline{\upxi_\epnu^{b,\gamma}},\overline{\upxi_\epnu^{c,\gamma}}\big)(z;z^-,z^+)=B_-\big(\dot{G}_\nu^{\overline{\gamma}},\upxi_\epnu^{\overline{b},\overline{\gamma}},\upxi_\epnu^{\overline{c},\overline{\gamma}}\big)(z;z^-,z^+)\,,
    \end{equs}
    and that similarly
    \begin{equs}
        \overline{B_-\big(\dot{G}_\nu^{\overline{\gamma}},\upxi_\epnu^{b,\gamma},\upxi_\epnu^{\overline{c},\gamma}\big)}(z;z^+,z^-)=B_+\big(\dot{G}_\nu^\gamma,\upxi_\epnu^{\overline{b},\overline{\gamma}},\upxi_\epnu^{c,\overline{\gamma}}\big)(z;z^-,z^+)\,.
    \end{equs}
    To conclude, it suffices to see that $a=b+c\Leftrightarrow\overline{a}=\overline{b}-\overline{c}$ and $a=b-c\Leftrightarrow\overline{a}=\overline{b}+\overline{c}$. The case of $(\upomega_\epmu^{a,\gamma})$ is dealt with identically. For $(\uptau_\epmu^{\a,\gamma})$ and $(\upzeta_\epmu^{\a,\gamma})$, the sketch of proof is the same. The only difference is that the induction given by the hierarchy is started with more base cases, because of the presence of the harmonic flow of the initial condition. We check only one example for conciseness.
    \begin{equs}
        \overline{\uptau_\epmu^{(0,1,1,0,1,0),\gamma}}(z;z^+,z^-)=\overline{\gamma}\lambda\delta(z-z_1^+)\delta(z-z_1^-)G_\mu^{\overline{\gamma}}\psi_\eps(z)=\uptau_\epmu^{(0,1,1,0,0,1),\overline{\gamma}}(z;z^-,z^+)\,.
    \end{equs}
    The same also holds for $(\upzeta_\epmu^{\a,\gamma})$, thus finishing the proof.
\end{proof}

\begin{lemma}\label{lemma:scalingzeta}
    The following holds for every $\a\in\mfM^\star_\pm$ and $\epmu$,
    \begin{equation}
        \upzeta_\epmu^{\a,\gamma}(\z)\eqlaw\gamma^{\underline{a}+1}\overline{\gamma}^{\overline{a}}\upzeta_\epmu^{\a,1}(\gamma\undcdot\z)\,,
    \end{equation}
    where we have set $\gamma\undcdot\z=(\gamma\cdot z;\gamma\cdot z^+,\overline\gamma\cdot z^-)$.
\end{lemma}
\begin{proof}
    The result is proved by induction on $\a$ with the help of the hierarchy \eqref{eq:hierarchyzeta}.\\
    \textit{Base case.} 
    The first base case is the noise, \ie $\a=(0,0,0,1,0,0)$, for which we have
    \begin{equs}
        \upzeta_\epmu^{(0,0,0,1,0,0),1}(\gamma z)&=\one_{\Re(\gamma)t>0}\sqrt{2}\xi_\eps(\gamma t,x)=\one_{t>0}\sqrt{2}\frac{\Re(\gamma)}{|\gamma|}\xi_\eps\big(\Re(\gamma)t,x\big)\eqlaw\one_{t>0}\frac{\sqrt{2\Re(\gamma)}}{|\gamma|}\xi_\eps(t,x)\\
        &\eqlaw\one_{t>0}\frac{\sqrt{2\Re(\gamma)}}{\gamma}\xi_\eps(t,x)\,.
    \end{equs}
    It then gives
    \begin{equs}
        \upzeta_\epmu^{(0,0,0,1,0,0),\gamma}(z)\eqlaw\gamma\upzeta_\epmu^{(0,0,0,1,0,0),1}(\gamma z)\,.
    \end{equs}
    The second base case is $\a=(0,2,1,0,0,0)$. In that case we have
    \begin{equs}
        \upzeta_\epmu^{(0,2,1,0,0,0),\gamma}(z;z_1^+,z_2^+,z_1^-)&=-\one_{t>0}\gamma^3\overline{\gamma}\lambda\times\delta\big(\gamma\cdot(z-z_1^+)\big)\delta\big(\gamma\cdot(z-z_2^+)\big)\delta\lambda\big(\overline\gamma\cdot(\overline{z}-z_1^-)\big)\\
        &=\gamma^3\overline{\gamma}\upzeta_\epmu^{(0,2,1,0,0,0),1}(\gamma\cdot z,\gamma\cdot z_1^+,\gamma\cdot z_2^+,\overline\gamma\cdot z_1^-)\,.
    \end{equs}
    Note that the choice of non-linearity that we made in equation \eqref{eq:eqholo} translates in a conjugation in time on $z$ in the third Dirac. To clear any doubt, we denote $\overline{z}=(\overline{\omega},x)$. All the remaining cases are the ones involving the harmonic flow of the initial condition. We treat only one of these terms and we let the reader fill in the details to get the others. We first need some basic fact on the harmonic flow of the initial condition. We have
    \begin{equs}
        M_\mu^\gamma\psi_\eps(z)&=\int_{\R^d}M_\mu^\gamma(t,x-y)u_0^\eps(y)\dint y=\int_{\R^d}M_\mu(\gamma t,x-y)u_0^\eps(y)\dint y\\
        &=M_\mu\psi_\eps(\gamma t,x)=M_\mu\psi_\eps(\gamma\cdot z)\,.
    \end{equs}
    It yields
    \begin{equs}
        \upzeta_\epmu^{(0,1,1,0,1,0),\gamma}(z;z_1^+,z_1^-)&=\gamma\lambda\times\delta(z-z_1^+)\delta(\overline{z}-z_1^-)G_\mu^\gamma\psi_\eps(z)\\
        &=\gamma^2\overline{\gamma}\lambda\times\delta\big(\gamma\cdot(z-z_1^+)\big)\delta\big(\overline{\gamma}\cdot(\overline{z}-z_1^-)\big)G_\mu\psi(\gamma\cdot z)\\
        &=\gamma^2\overline{\gamma}\upzeta_\epmu^{(0,1,1,0,1,0),1}(\gamma\cdot z,\gamma\cdot z_1^+,\overline \gamma\cdot z_1^-)\,.
    \end{equs}
    Finally, thanks to Proposition \ref{prop:scalingC}, we now that for every $n$, $\mfc_\eps^{3n-1,\gamma,\lambda}=\gamma\mfc_\eps^{3n-1,1,\lambda}$, thus treating all the terms that do not vanish at the boundary $\mu=0$.\\
    \textit{Induction.} We can now turn to the inductive part of the proof. Thanks to \eqref{eq:hierarchyzeta}, we can write 
    \begin{equs}
        \upzeta_\epmu^{\a,\gamma}=\upzeta_{\eps,0}^{\a,\gamma}-\int_0^\mu\sum_{\a=\b+\c} B_+\big(\dot{M}_\nu^\gamma,\upzeta_\epnu^{\b,\gamma},\upzeta_\epnu^{\c,\gamma}\big)+\sum_{\a=\b-\c}B_-\big(\dot{M}^{\overline{\gamma}}_\nu,\upzeta_\epnu^{\b,\gamma},\overline{\upzeta_\epnu^{\c,\gamma}}\big)\dint\nu\,.
    \end{equs}
    The term with $\mu=0$ is already tackled. We need to treat the two bilinear terms. Let us start with $B_+$. We have, for $\a,\b,\c$ such that $\a=\b+\c$,
    \begin{equs}
        ~&B_+\big(\dot M_\nu^\gamma,\upzeta_\epnu^{\b,\gamma},\upzeta_\epnu^{\c,\gamma}\big)(z;z_{1:\underline{a}}^+,z_{1:\overline{a}}^-)\\
        &=\sum_{j=1}^{\underline{a}}\int_{\Lambda^2}\upzeta_\epnu^{\b,\gamma}(z;z_{1:\underline{b}}^+,z_{1:\overline{b}}^-) \dot M_\nu^\gamma(z_j^+,w)\upzeta_\epnu^{\c,\gamma}(w;z_{\underline{b}+1:\underline{b}+\underline{c}}^+,z_{\overline{b}+1:\overline{b}+\overline{c}}^-)\dint z_j^+\dint w\\
        &\eqlaw\gamma^{\underline{b}+\underline{c}+2}\overline{\gamma}^{\overline{b}+\overline{c}}\int_{\Lambda^2}\upzeta_\epnu^{\b,1}(\gamma\cdot z;\gamma\cdot z_{1:\underline{b}}^+,\overline\gamma\cdot z_{1:\overline{b}}^-)\\
        &\qquad\qquad\qquad\qquad\qquad\qquad\qquad\quad\times\dot M_\nu(\gamma\cdot z_j^+,\gamma\cdot w)\upzeta_\epnu^{\c,1}(\gamma\cdot w;\gamma\cdot z_{\underline{b}+1:\underline{b}+\underline{c}}^+,\overline \gamma\cdot z_{\overline{b}+1:\overline{b}+\overline{c}}^-)\dint z_j^+\dint w\\
        &=\gamma^{\underline{b}+\underline{c}}\overline{\gamma}^{\overline{b}+\overline{c}}\int_{\Lambda^2}\upzeta_\epnu^{\b,1}(\gamma\cdot z;\gamma\cdot z_{1}^+,\dots,z_j^+,\dots,\gamma\cdot z_{\underline{b}}^+,\overline\gamma\cdot z_{1:\overline{b}}^-)\\
        &\qquad\qquad\qquad\qquad\qquad\qquad\qquad\qquad\qquad~~~\,\times\dot M_\nu( z_j^+, w)\upzeta_\epnu^{\c,1}(w;\gamma\cdot z_{\underline{b}+1:\underline{b}+\underline{c}}^+,\overline \gamma\cdot z_{\overline{b}+1:\overline{b}+\overline{c}}^-)\dint z_j^+\dint w\\
        &=\gamma^{\underline{a}+1}\overline{\gamma}^{\overline{a}}B_+\big(\dot M_\nu,\upzeta_\epnu^{\b,1},\upzeta_\epnu^{\c,1}\big)(\gamma\cdot z;\gamma\cdot z_{1:\underline{a}}^+,\overline \gamma\cdot z_{1:\overline{a}}^-)\,.
    \end{equs}
    For the $B_-$ term, take $\a,\b,\c$ such that $\a=\b-\c$. We use Lemma \ref{lemma:conj} to write $\overline{\upzeta_\epnu^{\c,\gamma}(z,z^+,z^-)}=\upzeta_\epnu^{\overline{c},\overline{\gamma}}(z,z^-,z^+)$. We get
    \begin{equs}
        ~&B_-\big(\dot M_\nu^{\overline{\gamma}},\upzeta_\epnu^{\b,\gamma},\overline{\upzeta_\epnu^{\c,\gamma}}\big)(z;z_{1:\underline{a}}^+,z_{1:\overline{a}}^-)\\
        &=\sum_{j=1}^{\overline{a}}\int_{\Lambda^2}\upzeta_\epnu^{\b,\gamma}(z;z_{1:\underline{b}}^+,z_{1:\overline{b}}^-) \dot M_\nu^{\overline{\gamma}}(z_j^-,w)\upzeta_\epnu^{\overline\c,\overline\gamma}(w;z_{\overline{b}+1:\overline{b}+\overline{c}}^-,z_{\underline{b}+1:\underline{b}+\underline{c}}^+)\dint z_j^-\dint w\\
        &\eqlaw\gamma^{\underline{b}+\overline{c}+1}\overline{\gamma}^{\overline{b}+\underline{c}+1}\int_{\Lambda^2}\upzeta_\epnu^{\b,1}(\gamma\cdot z;\gamma\cdot z_{1:\underline{b}}^+,\overline\gamma\cdot z_{1:\overline{b}}^-)\\
        &\qquad\qquad\qquad\qquad\qquad\qquad\qquad\quad\times\dot M_\nu(\overline\gamma\cdot z_j^+,\overline\gamma\cdot w)\upzeta_\epnu^{\overline\c,1}(\gamma\cdot w;\overline \gamma\cdot z_{\overline{b}+1:\overline{b}+\overline{c}}^-,\gamma\cdot z_{\underline{b}+1:\underline{b}+\underline{c}}^+)\dint z_j^-\dint w\\
        &=\gamma^{\underline{b}+\overline{c}+1}\overline{\gamma}^{\overline{b}+\underline{c}-1}\int_{\Lambda^2}\upzeta_\epnu^{\b,1}(\gamma\cdot z;\gamma\cdot z_{1:\underline{b}}^+,\overline\gamma\cdot z_{1}^-,\dots,z_j^-,\dots,\overline\gamma\cdot z_{\overline{a}}^-)\\
        &\qquad\qquad\qquad\qquad\qquad\qquad\qquad\qquad\qquad~~~\,\times\dot M_\nu(z_j^+, w)\upzeta_\epnu^{\overline\c,1}(w;\overline \gamma\cdot z_{\overline{b}+1:\overline{b}+\overline{c}}^-,\gamma\cdot z_{\underline{b}+1:\underline{b}+\underline{c}}^+)\dint z_j^-\dint w\\
        &=\gamma^{\underline{b}+\overline{c}+1}\overline{\gamma}^{\overline{b}+\underline{c}-1}\int_{\Lambda^2}\upzeta_\epnu^{\b,1}(\gamma\cdot z;\gamma\cdot z_{1:\underline{b}}^+,\overline\gamma\cdot z_{1}^-,\dots,z_j^-,\dots,\overline\gamma\cdot z_{\overline{a}}^-)\\
        &\qquad\qquad\qquad\qquad\qquad\qquad\qquad\qquad\qquad~~~\,\times\dot M_\nu(z_j^+, w)\overline{\upzeta_\epnu^{\c,1}}(w;\gamma\cdot z_{\underline{b}+1:\underline{b}+\underline{c}}^+,\overline \gamma\cdot z_{\overline{b}+1:\overline{b}+\overline{c}}^-)\dint z_j^-\dint w\\
        &=\gamma^{\underline{a}+1}\overline{\gamma}^{\overline{a}}B_-\big(\dot M_\nu,\upzeta_\epnu^{\b,1},\overline{\upzeta_\epnu^{\c,1}}\big)(\gamma\cdot z;\gamma\cdot z_{1:\underline{a}}^+,\overline \gamma\cdot z_{1:\overline{a}}^-)\,.
    \end{equs}
    This concludes the proof.
\end{proof}
We are eventually able to state and prove the main object of this section.
\begin{proposition}\label{prop:eqlawu}
    The following holds for every $(t,x)\in\Lambda$ and $\eps$
    \begin{equation}
        u^\gamma_\eps(t,x)\eqlaw u^1_\eps(\gamma t,x)\,,
    \end{equation}
    where $u_\eps^1$ is the solution to \eqref{eq:eqholobis} with $\gamma=1$.
\end{proposition}

\begin{proof}
    We now from \cite{Duc21,CF24a,FMM26} that the solution takes the following form
    \begin{equs}
        u_\eps^\gamma=M^\gamma\psi_\eps-\int_0^{T^{1/2}}\dot M_\nu^\gamma\big(F^\gamma_{\eps,T_\gamma^{1/2}}[0]+R^\gamma_{\eps,T_\gamma^{1/2}}\big)\dint\nu\,.
    \end{equs}
    Let us first take a look at the force term. Since the functional parameter is taken to be $0$, the only force coefficients that survive are the ones such that $\underline{a}=\overline{a}=0$. Lemma \ref{lemma:scalingzeta} directly gives $F_{\eps,T^{1/2}_\gamma}^\gamma[0](z)\eqlaw\gamma F_{\eps,T^{1/2}_\gamma}^1[0](\gamma\cdot z)$. If furthermore, we are able to show that $R_{\eps,T^{1/2}_\gamma}^\gamma(z)\eqlaw\gamma R_{\eps,T^{1/2}_\gamma}^1(\gamma\cdot z)$, we can conclude. In this last equality, there is an overload of notation, since the right-hand side is actually the remainder that we construct when solving the equation \eqref{eq:eqholochanged} in Section \ref{sec:system}. Indeed, we recall that we never really solve the SPDE \eqref{eq:eqholo} for complex times, as explained in the beginning of the section. We have also overloaded the notation for the stopping times, again using equality in law. We deliberately chose not to introduce further notation to preserve clarity. We conclude as follows.
    \begin{equs}
        \dot M_\nu^\gamma\big(F^\gamma_{\eps,T^{1/2}}[0]+R^\gamma_{\eps,T^{1/2}}\big)(t,x)&=\int_0^\infty\int_{\R^d}\dot M_\nu(\gamma t-\gamma s,x,y)\big(F^\gamma_{\eps,T_\gamma^{1/2}}[0]+R^\gamma_{\eps,T_\gamma^{1/2}}\big)(s,y)\dint y\dint s\\
        &\eqlaw\gamma\int_0^\infty\int_{\R^d}\dot M_\nu\big(\gamma t-\gamma s,x,y)(F^1_{\eps,T_\gamma^{1/2}}[0]+R^1_{\eps,T_\gamma^{1/2}}\big)(\gamma s,y)\dint y\dint s\\
        &=\int_0^\infty\int_{\R^d}\dot M_\nu(\gamma t-s,x,y)\big(F^1_{\eps,T_\gamma^{1/2}}[0]+R^1_{\eps,T_\gamma^{1/2}}\big)(s,y)\dint y\dint s\\
        &=\dot M_\nu\big(F^1_{\eps,T_\gamma^{1/2}}[0]+R^1_{\eps,T_\gamma^{1/2}}\big)(\gamma t,x)\,.
    \end{equs}
    We are left to prove the fact about the scaling of the remainder. For that purpose, the main observation is that, if $(v_\epmu^1,R_\epmu^1)$ solves the $+$ part of the system \eqref{eq:bigsys} with $\gamma=1$, then, the couple $\big(\tilde v_\epmu^\gamma(z),\tilde R_\epmu^\gamma(z)\big)\eqdef\big(v_\epmu^1(\gamma\cdot z),\gamma R_\epmu^1(\gamma\cdot z)\big)$ solves a posteriori the system
    \begin{subequations}
        \begin{empheq}[left=\empheqlbrace]{alignat=2}
          \D_\mu  \tilde v^{\gamma}_\epmu&=\dot M^{\gamma}_\mu\big(\tilde F^{\gamma}_{\eps,\mu}[\tilde v^{\gamma}_{\eps,\mu}]+\tilde R^{\gamma}_{\eps,\mu}\big)\,,\nonumber\\
          \D_\mu  \tilde R^{\gamma}_\epmu&=- \mathrm{D}\tilde F^{\gamma}_\epmu[\tilde v_\epmu^\gamma] \dot M^{\gamma}_\mu \tilde R^{\gamma}_\epmu-(\D_\mu+\mathrm{D}\tilde F_\epmu^\gamma[\tilde v_\epmu^\gamma]\dot M_\mu^\gamma)\tilde F_\epmu^\gamma[\tilde v^{\gamma}_{\eps,\mu}]\,,\nonumber
        \end{empheq}
    \end{subequations}
    where we have defined
    \begin{equs}
        \tilde F_\epmu^\gamma[\phi](z)=\sum_{\a\in \mfM^\star_+}\int_{\Lambda^{l(\a)}}\!\!\gamma^{\underline{a}+1}\overline{\gamma}^{\overline{a}}\upzeta_\epmu^{\a,1}(\gamma\cdot z;\gamma\cdot z^+,\overline{\gamma}\cdot z^-)\prod_{j=1}^{\underline{a}}\phi(z_j^+)\prod_{j=1}^{\overline{a}}\overline{\phi}(z_j^-)\dint z^+\dint z^-\,,
    \end{equs}
    so that $\tilde F_\epmu^\gamma\eqlaw F_\epmu^\gamma$ thanks to Lemma \ref{lemma:scalingzeta}. The main calculation to see the above is the following
    \begin{equs}
        \tilde F_\epmu^\gamma[\tilde v^{\gamma}_\epmu](z)&=\sum_{\a\in \mfM^\star_+}\int_{\Lambda^{l(\a)}}\!\!\gamma^{\underline{a}+1}\overline{\gamma}^{\overline{a}}\upzeta_\epmu^{\a,1}(\gamma\cdot z;\gamma\cdot z^+,\overline{\gamma}\cdot z^-)\prod_{j=1}^{\underline{a}}v_\epmu^1(\gamma\cdot z_j^+)\prod_{j=1}^{\overline{a}}\overline{v_\epmu^1}(\overline{\gamma}\cdot z_j^-)\dint z^+\dint z^-\\
        &=\gamma \sum_{\a\in \mfM^\star_+}\int_{\Lambda^{l(\a)}}\!\!\upzeta_\epmu^{\a,1}(\gamma\cdot z; z^+,z^-)\prod_{j=1}^{\underline{a}}v_\epmu^1(z_j^+)\prod_{j=1}^{\overline{a}}\overline{v_\epmu^1}(z_j^-)\dint z^+\dint z^-\\
        &=\gamma F_\epmu^1[v_\epmu^1](\gamma\cdot z)\,.
    \end{equs}
\end{proof}

\section{Cumulant analysis and renormalisation constants}\label{sec:proba}

\begin{definition}
Let $m\in\N$ and $X_1,\dots,X_m$ real-valued random variables. Their joint cumulant is defined as 
\begin{equation}
    \Kappa(X_1,\dots,X_m)=(-\i)^m\D_{t_1}\dots\D_{t_m}\log\left.\E\big[\exp(\i (t_1X_1+\dots+t_mX_m))\big]\right|_{t_1=\dots=t_m=0}.
\end{equation}
\end{definition}
Note that the latter is symmetric in its arguments, so that in only depends on the collection $\{X_i\colon i\in[m]\}$. We can now introduce the main object of this section.

\begin{definition}
    For $\mu\in [0,1]$, $\eps>0$, $m\in\N$, $a^m=(a_1,\dots,a_m)\in\mfM^m$, and $\gamma^m=(\gamma_1,\dots,\gamma_m)\in\{\gamma,\overline{\gamma}\}^m$, we let the cumulant
    \begin{equation}
    \Kappa_{\eps,\mu}^{a^m,\gamma^m}(\z_1,\dots,\z_m)=\Kappa\left(\upxi_{\eps,\mu}^{a_1,\gamma_1}(\z_1),\dots,\upxi_{\eps,\mu}^{a_m,\gamma_m}(\z_m)\right)\,.
    \end{equation}
    As cumulants are symmetric in their arguments, we point out that the order in $(a^m,\gamma^m)$ does not matter in any way.
\end{definition}

Before giving a norm for such objects, let us introduce some weights. We define $l(a^m)=\sum_{k=1}^m l(a_k)$. For $\w^m=(\w_1,\dots,\w_m)\in\Lambda^{l(a^m)+m}$, we set
\begin{equation}
    \langle\w^m\rangle_{\mu,m}=\prod_{k=1}^m\langle \w_k\rangle_\mu.
\end{equation}
We recall that we defined $\langle\bigcdot\rangle_\mu$ in \eqref{mubracket}.
\begin{definition}
    For $U\in\CS'(\Lambda^{l(a^m)+m})$, we set a norm, that is $L^\infty_{w_1}$ in $w_1$, $L^{1,p}_{t_i,x_i}$ in the $w_i$'s, where $w_i=(t_i,x_i)$, and $L^1_{t,w}$ in all the others ($w^{\pm}_i$). We have let for the rest of Section \ref{sec:proba}, $p=1\vee \frac{d}{2(\beta-1)}$.
    \begin{equation}
        \snnorm{U}_{\mu,m,\varsigma}=\sup_{w_1\in\Lambda}\int_{\R^{m-1}}\left(\int_{(\R^d)^{m-1}}\left(\int_{\Lambda^{m}}\langle\w^m\rangle_{\mu,m}^{\varsigma}\big|U(\w_1,\dots,\w_m)|\dint w^{+}_{1:m}\dint w^-_{1:m}\right)^{p}\dint x_{2:m}\right)^{\frac{1}{p}}\dint t_{2:m}\,.
    \end{equation}
\end{definition}
We also provide a notion of degree for cumulants, which reads
\begin{equation}
    \deg(a^m)=\sum_{k=1}^m\deg(a_k)+(m-1)\left(2+\frac{d}{p}\right)\,.
\end{equation}

\begin{definition}
    Let $\Lambda^{(l_+,l_-)}=\Lambda\times\Lambda^{l_+}\times\Lambda^{l_-}$ and let $U\in\mcS'\big(\prod_{i=1}^{m+1}\Lambda^{(l_+^i,l_-^i)}\big)$. For $1\leq k\leq m$, we define a distribution $A_{+}^{k}\big(G,U\big)\in\mcS'\big(\Lambda^{(l_+^k+l_+^{k+1}-1,l_-^k+l_-^{k+1})}\times\prod_{i\notin \{k,k+1\}}\Lambda^{(l_+^i,l_-^i)}\big)$ by
    \begin{equation}
        A_{+}^{k}\big(G,U\big)(\w_1,\dots,\w_m)=\sum_{j=1}^{l_+^k}\int_{\Lambda^2}U(\w_1,\dots,\w_k^{\downarrow j},\w_k^{\uparrow},\dots,\w_m)G(z_1,z_2)\dint z_1\dint z_2\,,
    \end{equation}
    where we have denoted $\w_k^{\downarrow j}=(w_k;w_{k,1}^+,\dots,z_1,\dots,w_{k,l_+^k-1}^+,w_{k,1:l_-^k}^-)$, with $z_1$ at position $j$, and $\w_k^{\uparrow}=(z_2;w_{k,l_+^k:l_+^k+l_+^{k+1}-1}^+,w_{k,l_-^k+1:l_-^k+l_-^{k+1}}^-)$.

    Similarly, for $m_1,m_2$ such that $m_1+m_2=m+1$, for $U\in\mcS'\big(\prod_{i=1}^{m_1+1}\Lambda^{(l_+^{1,i},l_-^{1,i})}\big)$ and $V\in\mcS'\big(\prod_{i=1}^{m_2+1}\Lambda^{(l_+^{2,i},l_-^{2,i})}\big)$, and for $1\leq k\leq m_1$, we define a distribution $B_+^k\big(G,U,V\big)\in\mcS'\big(\Lambda^{(l_+^{1,k}+l_+^{2,1}-1,l_-^{1,k}+l_-^{2,1})}\times\prod_{i\neq k}\Lambda^{(l_+^{1,i},l_-^{1,i})}\times\prod_{i\neq1}\Lambda^{(l_+^{2,i},l_-^{2,i})}\big)$ by
    \begin{multline}
        B_+^k\big(G,U,V\big)(\w_1,\dots,\w_m)=\sum_{j=1}^{l_+^{1,k}}\int_{\Lambda^2}U(\w_1,\dots,\w_k^{\downarrow j},\dots,\w_{m_1})\\\times G(z_1,z_2)V(\w_k^\uparrow,\w_{m_1+1},\dots,\w_{m_1+m_2+1})\dint z_1\dint z_2\,,
    \end{multline}
    where we have denoted $\w_k^{\downarrow j}=(w_k;w_{k,1}^+,\dots,z_1,\dots,w_{k,l_+^{1,k}-1}^+,w_{k,1:l_-^{1,k}}^-)$, with $z_1$ at position $j$, and $\w_k^{\uparrow}=(z_2;w_{k,l_+^{1,k}:l_+^{1,k}+l_+^{2,1}-1}^+,w_{k,l_-^{1,k}+1:l_-^{1,k}+l_-^{2,1}}^-)$.
    
    We let the reader fill in the details to get the definition of $A_-^k$ and $B_-^k$.
\end{definition}

It is now common in the flow approach literature that the hierarchy \eqref{eq:hierarchyxi}, combined with a right formula on cumulants of products, allows us to obtain the following hierarchy on the cumulants of the cumulants of the stationary force coefficients.

\begin{lemma}\label{lemma:decompderivative}
    $\D_\mu\Kappa_\epmu^{\a^m,\gamma^m}$ is a linear combination of terms of the form
    \begin{equation}
        A_\pm^k\big(\dot G_\mu^\gamma,\Kappa_\epmu^{b^{m+1},\gamma^{m+1}}\big)\,,~\mathrm{and}~~B_\pm^k\big(\dot G_\mu^\gamma,\Kappa_\epmu^{b^{m_1},\gamma^{m_1}},\Kappa_\epmu^{c^{m_2},\gamma^{m_2}}\big)\,,
    \end{equation}
    where
    \begin{enumerate}
        \item For $A_+$, $\gamma^{m+1}=(\gamma,\gamma^m)$ and
        \begin{equation}
            \sum_{k=1}^m a_{\circ,k}=\sum_{k=1}^{m+1}b_{\circ,k}+1,\quad\sum_{k=1}^m \underline{a}_k=\sum_{k=1}^{m+1}\underline{b}_k-1,\quad\sum_{k=1}^m \overline{a}_k=\sum_{k=1}^{m+1}\overline{b}_k\,.
        \end{equation}
        \item For $B_+$, $m_1+m_2=m+1$ and $\gamma^{m_1}$ and $\gamma^{m_2}$ are such that $\gamma^m=(\gamma, \gamma^{m_1}, \gamma^{m_2})$, and
        \begin{equation}
            \sum_{k=1}^m a_{\circ,k}=\sum_{k=1}^{m_1} b_{\circ,k}+\sum_{k=1}^{m_2} c_{\circ,k}+1,\quad\sum_{k=1}^m \underline{a}_k=\sum_{k=1}^{m_1} \underline{b}_k+\sum_{k=1}^{m_2} \underline{c}_k-1,\quad\sum_{k=1}^m \overline{a}_k=\sum_{k=1}^{m_1} \overline{b}_k+\sum_{k=1}^{m_2} \overline{c}_k\,,
        \end{equation}
        \item For $A_-$, $\gamma^{m+1}=(\overline\gamma,\gamma^m)$ and
        \begin{equation}
            \sum_{k=1}^m a_{\circ,k}=\sum_{k=1}^{m+1}b_{\circ,k}+1,\quad\sum_{k=1}^m \underline{a}_k=\sum_{k=1}^{m+1}\underline{b}_k,\quad\sum_{k=1}^m \overline{a}_k=\sum_{k=1}^{m+1}\overline{b}_k-1\,,
        \end{equation}
        \item For $B_-$, $m_1+m_2=m+1$ and $\gamma^{m_1}$ and $\gamma^{m_2}$ are such that $\gamma^m=(\overline\gamma, \gamma^{m_1}, \gamma^{m_2})$, and
        \begin{equation}
            \sum_{k=1}^m a_{\circ,k}=\sum_{k=1}^{m_1} b_{\circ,k}+\sum_{k=1}^{m_2} c_{\circ,k}+1,\quad\sum_{k=1}^m \underline{a}_k=\sum_{k=1}^{m_1} \underline{b}_k+\sum_{k=1}^{m_2} \underline{c}_k,\quad\sum_{k=1}^m \overline{a}_k=\sum_{k=1}^{m_1} \overline{b}_k+\sum_{k=1}^{m_2} \overline{c}_k-1\,,
        \end{equation}
    \end{enumerate}
\end{lemma}

This decomposition is used to show by induction the following bound on the cumulants, using Lemma \ref{estimdotGmu}
\begin{equation}\label{eq:boundcumu}
    \nnorm{K_\mu^{\otimes(l(a^m)+m)}\Kappa_\epmu^{a^m,\gamma^m}}_{\mu,m,\varsigma}\lesssim\mu^{\deg(a^m)}\,,
\end{equation}
for any $\varsigma>0$. We do not reprove this fact in the present work since it is now quite classical. We refer the reader to \cite{DGR23} or \cite{FMM26}. However, we want to point out some facts that are proper to this model. The first one is that the choice of norm slightly differs from the references we gave. This is due to the covariance structure of the noise. The choice of norm is here closer to the one made in \cite{CF24a}, in which the noise is also white in time. It is not hard to see that the proof can be adapted to this choice. The other important fact is about which terms need to be renormalised. We recall that we renormalise only the \textit{relevant} cumulants, \ie the ones of negative degree. The first fact is that $p$ is chosen such that
\begin{equs}
    \deg(a^m)\geq -m\beta+2(m-1)\beta=(m-2)\beta\,,
\end{equs}
so $\deg(a^m)<0$ imposes $m=1$. Then, we check that, for $a\in\mfM_+$,
\begin{equs}
    \deg(a)<0\Rightarrow \underline{a}+\overline{a}\in\{0,1\}\,.
\end{equs}
Since $\E\Big[\prod_{i=1}^{l_+}\xi(z_i^+)\prod_{j=1}^{l_-}\overline{\xi(z_j^-)}\Big]$ is non-zero if and only if $l_+=l_-$, we also get that $a_\pm=0$ which in turn implies that $\overline{a}-\underline{a}+1=0$. This and the fact that in the above the case of $0$ also give a null expectation (as the expectation of an odd product of Gaussian random variables vanishes), we get $\underline{a}=1$ and $\overline{a}=0$. It gives the form of the renormalisation of the SPDE \eqref{eq:bareeq}. One can deduce a similar statement for $a\in\mfM_-$. We can now show that the cumulants and the renormalisation constants satisfy a nice scaling property. 

\begin{proposition}\label{prop:scalingC}
    Let $m\in\N$ and fix some $a^m$ and $\gamma^m$. For $\z=(z,z^+,z^-)\in\Lambda\times\Lambda^{l_+}\times\Lambda^{l_-}$, we define $\gamma\undcdot\z=(\gamma\cdot z,\gamma\cdot z^+,\overline\gamma\cdot z^-)$ and $\gamma\ovcdot\z=(\overline\gamma\cdot z,\gamma\cdot z^+,\overline\gamma\cdot z^-)$, where we recall that, for $z=(t,x)$, $\gamma\cdot z=(\gamma t,x)$. Moreover, we set $\bigcdot_i=~\undcdot$ if $\gamma_i=\gamma$ and $\bigcdot_i=\ovcdot$ if $\gamma_i=\overline\gamma$. Then, the following holds for every $a^m$,
    \begin{equation}\label{eq:scalingcumu}
        \Kappa_\epmu^{a^m,\gamma^m}(\z_1,\dots,\z_m)=\prod_{k=1}^m\gamma^{\underline{a}_k}\overline{\gamma}^{\overline{a}_k}\gamma_k\times\Kappa_\epmu^{a^m,1^m}(\gamma_1\bigcdot_1\z_1,\dots,\gamma_m\bigcdot_m\z_m)\,.
    \end{equation}
    We have denoted $1^m=(1,\dots,1)\in\C^m$.
    Moreover,
    \begin{equation}\label{eq:scalingC}
        \mfC_\eps^{\gamma,\lambda}=\gamma\mfC_\eps^{1,\lambda}\,.
    \end{equation}
\end{proposition}

\begin{proof}
    More precisely, we show that for every $n$, $\mfc_\eps^{3n-1,\gamma,\lambda}=\gamma\mfc_\eps^{3n-1,1,\lambda}$. For that, we are going to show that the following slight generalisation of \eqref{eq:scalingcumu} holds for $s\in\{0,1\}$
    \begin{equs}\label{eq:HRcum}
        \D_\mu^s\Kappa_\epmu^{a^m,\gamma^m}(\z_1,\dots,\z_m)=\prod_{k=1}^m\gamma^{\underline{a}_k}\overline{\gamma}^{\overline{a}_k}\gamma_k\times\D_\mu^s\Kappa_\epmu^{a^m,1^m}(\gamma_1\bigcdot_1\z_1,\dots,\gamma_m\bigcdot_m\z_m)\,.
    \end{equs}
    The proof of the statement on the renormalisation constants is going to be intertwined with the one of this latter fact on the cumulants. The proof works by induction.\\
    \textit{Base case.} The first base case is $a^1=(0,2,1,0)$, for which
    \begin{equs}
        \Kappa_\epmu^{(0,2,1,0),\gamma}(z,z_1^+,z_2^+,z_1^-)=-\gamma\lambda\times\delta(z-z_1^+)\delta(z-z_2^+)\delta(\overline{z}-z_1^-)\,.
    \end{equs}
    Using the fact that $\delta(\gamma\cdot z)=\gamma^{-1}\delta(z)$, we get that
    \begin{equs}
        \Kappa_\epmu^{(0,2,1,0),\gamma}(z,z_1^+,z_2^+,z_1^-)&=-\gamma^3\overline{\gamma}\lambda\times\delta\big(\gamma\cdot(z-z_1^+)\big)\delta\big(\gamma\cdot(z-z_2^+)\big)\delta\big(\overline\gamma\cdot(\overline{z}-z_1^-)\big)\\
        &=\gamma^3\overline{\gamma}\times\Kappa_\epmu^{(0,2,1,0),1}(\gamma\cdot z,\gamma\cdot z_1^+,\gamma\cdot z_2^+,\overline\gamma\cdot z_1^-)\,.
    \end{equs}
    With the same computation, we also get that
    \begin{equs}
        \Kappa_\epmu^{(0,1,2,0),\overline\gamma}(z,z_1^+,z_1^-,z_2^-)&=\overline{\gamma}^3\gamma\times\Kappa_\epmu^{(0,1,2,0),1}(\overline\gamma\cdot z,\gamma\cdot z_1^+,\overline\gamma\cdot z_1^-,\overline\gamma\cdot z_2^-)\,.
    \end{equs}
    The second base case is $a^2=\big((0,0,0,1),(0,0,0,-1)\big)$, for which
    \begin{equs}
        \Kappa_\epmu^{((0,0,0,1),(0,0,0,-1)),\gamma}(z_1,z_2)=2\Re(\gamma)\E[\xi_\eps(z_1)\overline{\xi_\eps(z_2)}]\,.
    \end{equs}
    Thanks to the extension of the noise to complex times that we defined in \eqref{eq:extensiondistribs}, we have, for $z=(t,x)$,
    \begin{equs}
        \xi_\eps(\gamma\cdot z)=\frac{\Re(\gamma)}{|\gamma|}\xi_\eps(\Re(\gamma)t,x)\,,~\mathrm{and}~~\overline{\xi_\eps(\overline\gamma\cdot z)}=\frac{\Re(\gamma)}{|\gamma|}\overline{\xi_\eps(\Re(\gamma)t,x)}\,.
    \end{equs}
    It eventually yields
    \begin{equs}
        \gamma\overline{\gamma}\Kappa_\epmu^{((0,0,0,1),(0,0,0,-1)),(1,1)}(\gamma\cdot z_1,\overline\gamma\cdot z_2)&=2\Re(\gamma)^2\E\big[\xi_\eps(\Re(\gamma)t_1,x_1)\overline{\xi_\eps(\Re(\gamma)t_2,x_2)}\big]\\
        &=2\Re(\gamma)^2\delta\big(\Re(\gamma)(t_1-t_2)\big)\mathrm{Cov}_\eps(x_1-x_2)\\
        &=2\Re(\gamma)\delta(t_1-t_2)\mathrm{Cov}_\eps(x_1-x_2)\\
        &=2\Re(\gamma)\E\big[\xi_\eps(t_1,x_1)\overline{\xi_\eps(t_2,x_2)}\big]\\
        &=\Kappa_\epmu^{((0,0,0,1),(0,0,0,-1)),(\gamma,\overline{\gamma})}(z_1,z_2)\,.
    \end{equs}
    We have denoted $\mathrm{Cov}_\eps(x)=(\rho_\eps^{*2}*\mathrm{Cov})(x)$, with $\rho_\eps(x)=\eps^{-d}\rho(x/\eps)$. We can clearly see here why we took a mollification in space only.\\
    \textit{Induction.} We now turn to the induction step. We write
    \begin{equs}\label{eq:decompKmu}
        \Kappa_\epmu^{a^m,\gamma}(\z^m)=\Kappa_{\eps,0}^{a^m,\gamma}(\z^m)+\int_0^\mu \D_\nu \Kappa_\epnu^{a^m,\gamma}(\z^m)\dint \nu\,.
    \end{equs}
    The first term is non-zero only in the base case, that is already dealt with, and for $a^1=(3n-1,1,0,0)$, in which case
    \begin{equs}
        \Kappa_{\eps,0}^{(3n-1,1,0,0),\gamma}(z,z^+)=\mfc_\eps^{n,\gamma,\lambda}\delta(z-z^+)\,.
    \end{equs}
    We know from \cite{DGR23} that $\mfc_\eps^{3n-1,\gamma}$ has a rather explicit form given by
    \begin{equs}
        \mfc_\eps^{3n-1,\gamma,\lambda}=\int_0^1 \I\big(K_\nu^{\otimes2}\D_\nu\Kappa_{\eps,\nu}^{(3n-1,1,0,0),\gamma}\big)\dint\nu\,.
    \end{equs}
    $\I$ is a localisation operator given, for $U\in\mcS'(\Lambda^2)$, by
    \begin{equs}
        \I U(z_1)=\int_\Lambda U(z_1,z_2)\dint z_2\,.
    \end{equs}
    A direct computation gives
    \begin{equs}
        \mfc_\eps^{3n-1,\gamma,\lambda}=\int_\Lambda \D_\nu\Kappa_{\eps,\nu}^{(3n-1,1,0,0),\gamma}(0,z^+)\dint z^+\,.
    \end{equs}
    We can now use our induction hypothesis to get
    \begin{equs}
        \mfc_\eps^{3n-1,\gamma,\lambda}=\gamma^2\int_\Lambda \D_\nu\Kappa_{\eps,\nu}^{(3n-1,1,0,0),1}(0,\gamma\cdot z^+)\dint z^+=\gamma\int_\Lambda \D_\nu\Kappa_{\eps,\nu}^{(3n-1,1,0,0),1}(0,z^+)\dint z^+=\gamma \mfc_\eps^{3n-1,1,\lambda}\,.
    \end{equs}
    Finally, using the scaling properties of the Dirac $\delta$, we have
    \begin{equs}
         \Kappa_{\eps,0}^{(3n-1,1,0,0),\gamma}(z,z^+)=\gamma^2\Kappa_{\eps,0}^{(3n-1,1,0,0),1}(\gamma\cdot z,\gamma\cdot z^+)\,.
    \end{equs}
    To finish the proof, we need to treat the derivative term in \eqref{eq:decompKmu}. For that, we need to use Lemma \ref{lemma:decompderivative}. Thanks to the induction hypothesis \eqref{eq:HRcum} and computations on $A_\pm$ and $B_\pm$ that are in every way similar to what we did with $B_\pm$ in the proof of Lemma \ref{lemma:scalingzeta}, one can show the four following facts, using the notations of Lemma \ref{lemma:decompderivative}.
    \begin{enumerate}
        \item For $A_+$, we have
        \begin{equs}
            ~&A_+^k\big(\dot G_\mu^\gamma,\Kappa_\epmu^{b^{m+1},\gamma^m}\big)(\z_1,\dots,\z_m)\\
            &=\gamma^{-2}\gamma^{\sum_{k=1}^{m+1}\underline{b}_k}\overline{\gamma}^{\sum_{k=1}^{m+1}\overline{b}_k}\prod_{k=1}^{m+1}\gamma_k\times A_+^k\big(\dot G_\mu,\Kappa_\epmu^{b^{m+1},1^{m+1}}\big)(\gamma_1\bigcdot_1\z_1,\dots,\gamma_m\bigcdot_m\z_m)\\
            &=\gamma^{\sum_{k=1}^{m}\underline{a}_k}\overline{\gamma}^{\sum_{k=1}^{m+1}\overline{a}_k}\prod_{k=1}^m\gamma_k\times A_+^k\big(\dot G_\mu,\Kappa_\epmu^{b^{m+1},1^{m+1}}\big)(\gamma_1\bigcdot_1\z_1,\dots,\gamma_m\bigcdot_m\z_m)\,.
        \end{equs}
        Note that we have made a slight abuse of notation by not highlighting the dependence of $\gamma^m$ and $\gamma^{m+1}$ on $a^m$ and $b^{m+1}$ respectively. Also note that a factor $\gamma$ gets out of the product in that case.
        \item For $A_-$, we have 
        \begin{equs}
            ~&A_-^k\big(\dot G_\mu^\gamma,\Kappa_\epmu^{b^{m+1},\gamma^m}\big)(\z_1,\dots,\z_m)\\
            &=\overline{\gamma}^{-2}\gamma^{\sum_{k=1}^{m+1}\underline{b}_k}\overline{\gamma}^{\sum_{k=1}^{m+1}\overline{b}_k}\prod_{k=1}^{m+1}\gamma_k\times A_-^k\big(\dot G_\mu,\Kappa_\epmu^{b^{m+1},1^{m+1}}\big)(\gamma_1\bigcdot_1\z_1,\dots,\gamma_m\bigcdot_m\z_m)\\
            &=\gamma^{\sum_{k=1}^{m}\underline{a}_k}\overline{\gamma}^{\sum_{k=1}^{m}\overline{a}_k}\prod_{k=1}^m\gamma_k\times A_-^k\big(\dot G_\mu,\Kappa_\epmu^{b^{m+1},1^{m+1}}\big)(\gamma_1\bigcdot_1\z_1,\dots,\gamma_m\bigcdot_m\z_m)\,.
        \end{equs}
        This time it is a factor $\overline{\gamma}$ that gets out of the product.
        \item For $B_+$ we have
        \begin{equs}
            ~&B_+^k\big(\dot G_\mu^\gamma,\Kappa_\epmu^{b^{m_1},\gamma^{m_1}},\Kappa_\epmu^{c^{m_2},\gamma^{m_2}}\big)(\z_1,\dots,\z_m)\\
            &=\gamma^{-2}\gamma^{\sum_{k=1}^{m_1}\underline{b}_k+\sum_{k=1}^{m_2}\underline{c}_k}\overline{\gamma}^{\sum_{k=1}^{m_1}\overline{b}_k+\sum_{k=1}^{m_2}\overline{c}_k}\prod_{k=1}^{m_1}\gamma_k\prod_{k=1}^{m_2}\gamma_k\\
            &\qquad \times B_+^k\big(\dot G_\mu,\Kappa_\epmu^{b^{m_1},1^{m_1}},\Kappa_\epmu^{c^{m_2},1^{m_2}}\big)(\gamma_1\bigcdot_1\z_1,\dots,\gamma_m\bigcdot_m\z_m)\\
            &=\gamma^{\sum_{k=1}^{m}\underline{a}_k}\overline{\gamma}^{\sum_{k=1}^{m}\overline{a}_k}\prod_{k=1}^m\gamma_k\times B_+^k\big(\dot G_\mu,\Kappa_\epmu^{b^{m_1},1^{m_1}},\Kappa_\epmu^{c^{m_2},1^{m_2}}\big)(\gamma_1\bigcdot_1\z_1,\dots,\gamma_m\bigcdot_m\z_m)\,.
        \end{equs}
        We have made the same abuse of notation for $\gamma^m$.
        \item Finally, for $B_-$, we have
        \begin{equs}
            ~&B_-^k\big(\dot G_\mu^\gamma,\Kappa_\epmu^{b^{m_1},\gamma^{m_1}},\Kappa_\epmu^{c^{m_2},\gamma^{m_2}}\big)(\z_1,\dots,\z_m)\\
            &=\overline{\gamma}^{-2}\gamma^{\sum_{k=1}^{m_1}\underline{b}_k+\sum_{k=1}^{m_2}\underline{c}_k}\overline{\gamma}^{\sum_{k=1}^{m_1}\overline{b}_k+\sum_{k=1}^{m_2}\overline{c}_k}\prod_{k=1}^{m_1}\gamma_k\prod_{k=1}^{m_2}\gamma_k\\
            &\qquad \times B_-^k\big(\dot G_\mu,\Kappa_\epmu^{b^{m_1},1^{m_1}},\Kappa_\epmu^{c^{m_2},1^{m_2}}\big)(\gamma_1\bigcdot_1\z_1,\dots,\gamma_m\bigcdot_m\z_m)\\
            &=\gamma^{\sum_{k=1}^{m}\underline{a}_k}\overline{\gamma}^{\sum_{k=1}^{m}\overline{a}_k}\prod_{k=1}^m\gamma_k\times B_-^k\big(\dot G_\mu,\Kappa_\epmu^{b^{m_1},1^{m_1}},\Kappa_\epmu^{c^{m_2},1^{m_2}}\big)(\gamma_1\bigcdot_1\z_1,\dots,\gamma_m\bigcdot_m\z_m)\,.
        \end{equs}
    \end{enumerate}
    Putting everything together, we get
    \begin{equs}
        \D_\mu\Kappa^{a^m,\gamma^m}_\epmu(\z_1,\dots,\z_m)=\prod_{k=1}^m\gamma^{\underline{a}_k}\overline{\gamma}^{\overline{a}_k}\gamma_k\times\D_\mu\Kappa_\epmu^{a^m,1^m}(\gamma_1\bigcdot_1\z_1,\dots,\gamma_m\bigcdot_m\z_m)\,,
    \end{equs}
    which allows to conclude.
\end{proof}

\begin{appendix}

\renewcommand*{\thetheorem}{\mbox{\textrm A.\arabic{theorem}}}

\section{Estimates on the kernels}

The proofs of the following lemmas can be found in \cite[Appendix A]{FMM26}.

\begin{lemma}\label{estimdotGmu}
The following holds uniformly in $\mu\in(0,1]$, $\varsigma\geq0$, and $p,q\in[1,\infty]$,
    \begin{equation}
        \big\|\langle\bigcdot\rangle_\mu^{\varsigma}\mcR_\mu\dot G_\mu^\gamma\big\|_{L_{t,x}^{p,q}}\lesssim\mu^{-1+2/p-d(1-1/q)}\,,
    \end{equation}
\end{lemma}

\begin{lemma}\label{lemma:estimMmu}
    For every $n,m$, $\varsigma>0$, and $\mu\in(0,1]$,
    \begin{equation}
        \big\|\langle\bigcdot\rangle^\varsigma_\mu(\mcT_{\mu,n}\otimes \mcT_{\mu,m})\dot M_\mu^\gamma\big\|_{\mcL(L_{t,x}^{\infty,p})}\lesssim\mu\,,
    \end{equation}
    and
    \begin{equation}\label{eq:PmuTmuMmu}
        \big\|(\msP_\mu\otimes \mcT_{\mu,m})\dot M_\mu^\gamma\big\|_{\mcL(L_{t,x}^{\infty,p})}\lesssim\mu\,,
    \end{equation}
    $\mcT_{\mu,n}$ was defined in \eqref{eq:defTmu}.
\end{lemma}

\begin{lemma}\label{lemma:diffMG}
    The following holds for $\mu\in(0,1]$
    \begin{equation}
        \big|\mcT_{\mu,0}^{\otimes2}\big(\dot M_\mu^\gamma-\mcb W_\mu^{-1}\dot G_\mu^\gamma\big)(t,x,y)\big|\lesssim \mu^2\big(\mu^2+|x|^2+|y|^2\big)|\dot G_\mu^\gamma(8t,x-y)|\,.
    \end{equation}
\end{lemma}

\section{Mild formulation of the problem in complex times}\label{sec:system}

This section is devoted to the resolution of \eqref{eq:eqholobis} in complex times, in the sense of the system \eqref{eq:eqholochanged}. Let us explain first how to get the stochastic objects. Thanks to Lemma \ref{lemma:conj}, the ansatz for \eqref{eq:eqholoray} is
\begin{align}
    F^{+;\gamma}_\epmu[u_+,u_-](z)&=\sum_{\a\in \mfM^\star_+}\int_{\Lambda^{l(\a)}}\!\!\tilde\upzeta_\epmu^{\a,\gamma}(z;z_{1:\underline{a}}^+,z_{1:\overline{a}}^-)\prod_{j=1}^{\underline{a}}u_+(z_j^+)\prod_{j=1}^{\overline{a}}u_-(z_j^-)\dint z^+_{1:\underline{a}}\dint z^-_{1:\overline{a}}\,,\\
    F^{-;\overline\gamma}_\epmu[u_+,u_-](z)&=\sum_{\a\in \mfM^\star_-}\int_{\Lambda^{l(\a)}}\!\!\tilde\upzeta_\epmu^{\a,\overline\gamma}(z;z_{1:\underline{a}}^+,z_{1:\overline{a}}^-)\prod_{j=1}^{\underline{a}}u_-(z_j^+)\prod_{j=1}^{\overline{a}}u_+(z_j^-)\dint z^+_{1:\underline{a}}\dint z^-_{1:\overline{a}}\,,
\end{align}
where the coefficients $\tilde\upzeta_\epmu^{\a,\gamma}$ are defined like the effective force coefficients $\upzeta_\epmu^{\a,\gamma}$ with the change of initial condition $\tilde\upzeta_\epmu^{(0,0,0,1),\gamma}(z)=\gamma\xi_\eps(\gamma \cdot z)$.
Note that, thanks to Lemma \ref{lemma:scalingzeta}, $\tilde\upzeta_\epmu^{\a,\gamma}\eqlaw\upzeta_\epmu^{\a,\gamma}$. Such force coefficients are already constructed, following Proposition \ref{prop:estimatezeta}.
The flow system that we consider to solve \eqref{eq:eqholochanged} is the following, for which one has to think that $\gamma$ is in fact $\gamma\e^{\i\theta}$.
    \begin{empheq}[left=\empheqlbrace]{equation}\label{eq:bigsys}
    \begin{aligned}
        \D_\mu  v^{+;\gamma}_\epmu&=\dot M^{\gamma}_\mu\big(F^{+;\gamma}_{\eps,\mu}[v^{+;\gamma}_{\eps,\mu},v^{-;\gamma}_{\eps,\mu}]+R^{+;\gamma}_{\eps,\mu}\big)\,,\\
        \D_\mu  v^{-;\gamma}_\epmu&=\dot M^{\overline\gamma}_\mu\big(F^{-;\overline\gamma}_{\eps,\mu}[v^{+;\gamma}_{\eps,\mu},v^{-;\gamma}_{\eps,\mu}]+R^{-;\gamma}_{\eps,\mu}\big)\,,\\
        \D_\mu  R^{+;\gamma}_\epmu&=- \mathrm{D}F^{+;\gamma}_\epmu[v^{+;\gamma}_{\eps,\mu},v^{-;\gamma}_{\eps,\mu}] \dot M^{\gamma}_\mu R^{+;\gamma}_\epmu-(\D_\mu+\mathrm{D}F^{+;\gamma}_{\eps,\mu}[v^{+;\gamma}_{\eps,\mu},v^{-;\gamma}_{\eps,\mu}]\dot M_\mu^\gamma)F^{+;\gamma}_{\eps,\mu}[v^{+;\gamma}_{\eps,\mu},v^{-;\gamma}_{\eps,\mu}]\,\\
        \D_\mu  R^{-;\gamma}_\epmu&=- \mathrm{D}F^{-;\overline\gamma}_\epmu[v^{+;\gamma}_{\eps,\mu},v^{-;\gamma}_{\eps,\mu}] \dot M^{\overline\gamma}_\mu R^{-;\gamma}_\epmu-(\D_\mu+\mathrm{D}F^{-;\overline\gamma}_{\eps,\mu}[v^{+;\gamma}_{\eps,\mu},v^{-;\gamma}_{\eps,\mu}]\dot M_\mu^{\overline\gamma})F^{-;\overline\gamma}_{\eps,\mu}[v^{+;\gamma}_{\eps,\mu},v^{-;\gamma}_{\eps,\mu}]\,.
    \end{aligned}
    \end{empheq}
The resolution of such a system is done with the help of a fixed-point argument, which is now quite classical in the flow approach literature (see \cite{CF24a,FMM26}), and we do not reprove it. The result is the lemma below.
\begin{lemma}
    Let $\gamma\in\C_+$. There is a random $T^{\gamma}$ such that the map implicitly defined by the system above is a contraction for the norm
    \begin{multline}
        \snnorm{v^+,v^-,R^+,R^-}=\|u_0^\eps\|_{\msB_{p,q}^{-\sigma}}^{-1}\sup_{\mu\in (0,\sqrt{T^{\gamma}}]}\sup_{t\leq T^{\gamma}}\Big(\|\mcT_\mu v_\mu^+(t)\|_{L^p_x}\vee\|\mcT_\mu v_\mu^-(t)\|_{L^p_x}\\
        \vee\|J_\mu R_\mu^+(t)\|_{L^p_x}\vee\|J_\mu R_\mu^-(t)\|_{L^p_x}\Big)\,,
    \end{multline}
    and thus admits a unique fixed point $(v_\epmu^{+;\gamma},v_\epmu^{-;\overline{\gamma}},R_\epmu^{+;\gamma},R_\epmu^{-;\overline{\gamma}})$.
\end{lemma}
The uniqueness also yields the relations
\begin{equs}
    v_\epmu^{-;\overline{\gamma}}(\omega,x)=\overline{v_\epmu^{+;\gamma}(\overline{\omega},x)}\,,~\mathrm{and}~~R_\epmu^{+;\gamma}(\omega,x)=\overline{R_\epmu^{-;\overline{\gamma}}(\overline{\omega},x)}\,.
\end{equs}
The solution is then constructed as so
\begin{empheq}[left=\empheqlbrace]{equation}
    \begin{aligned}
    u_+^{\gamma\e^{\i\theta}}(z)&=M^{\gamma\e^{\i\theta}}\psi_\eps-\int_0^{T_{\gamma\e^{\i\theta}}^{1/2}}\dot M_\nu^{\gamma\e^{\i\theta}}(F^{+;\gamma\e^{\i\theta}}_{\eps,T_{\gamma\e^{\i\theta}}^{1/2}}[0]+R^{+;\gamma\e^{\i\theta}}_{\eps,T_{\gamma\e^{\i\theta}}^{1/2}})\dint\nu\,,\\
    u_-^{\overline\gamma\e^{-\i\theta}}(z)&=M^{\overline\gamma \e^{-\i\theta}}\overline{\psi_\eps}-\int_0^{T_{\gamma\e^{\i\theta}}^{1/2}}\dot M_\nu^{\overline\gamma \e^{-\i\theta}}(F^{-;\gamma\e^{-\i\theta}}_{\eps,T_{\gamma\e^{\i\theta}}^{1/2}}[0]+R^{-;\overline\gamma\e^{-\i\theta}}_{\eps,T_{\gamma\e^{\i\theta}}^{1/2}})\dint\nu\,.
    \end{aligned}
\end{empheq}
We have again thanks to uniqueness the relation $u_-^{\overline\gamma\e^{-\i\theta}}(\omega,x)=\overline{u_+^{\gamma\e^{\i\theta}}(\overline{\omega},x)}$.

\end{appendix}

\begingroup
\tiny
\bibliographystyle{alpha}
\bibliography{BiblioGL}
\endgroup

\end{document}